\documentclass[journal,twoside,web]{ieeecolor}
\usepackage{xcolor}
\usepackage{generic}
\usepackage[centertags]{amsmath}
\usepackage{tikz}
\usepackage{amsfonts}
\usepackage{newlfont}
\usepackage{graphicx}
\usepackage{epstopdf}
\usepackage[utf8]{inputenc}
\usepackage{indentfirst}
\usepackage{float}
\usepackage{cite}
\usepackage{CJK}
\usepackage{array}
\usepackage{amsfonts}
\usepackage{mathrsfs}
\usepackage{amssymb}
\usepackage{mathtools}

\usepackage{cancel}
\usepackage[normalem]{ulem}
\usepackage{tcolorbox}

\usepackage{booktabs}
\usepackage{tabularx}
\usepackage{array}

\usepackage{url}
\usepackage{hyperref}

\newtheorem{thm}{Theorem}[section]

\newtheorem{lem}[thm]{Lemma}

\newtheorem{defn}[thm]{Definition}
\newtheorem{assum}[thm]{Assumption}

\newtheorem{rem}[thm]{Remark}

\def\norm#1{\|{#1}\|}

\newcommand{\Reg}{\mathrm{Reg}}
\newcommand{\Exp}{\operatorname{Exp}}
\newcommand{\Log}{\operatorname{Log}}
\newcommand{\grad}{\operatorname{grad}}
\newcommand{\esgrad}{\mathbf{g}^{\delta}}
\newcommand{\ptg}[2]{\Gamma^{#2}_{#1}}
\newcommand{\kmin}{K_{\min}}
\newcommand{\kmax}{K_{\max}}

\newcommand{\sumstep}{\sum_{t=1}^T\eta_t}

\newcolumntype{C}[1]{>{\centering\arraybackslash}p{#1}}

\allowdisplaybreaks

\begin{document}

\title{Decentralized Online Riemannian Optimization for Strongly Geodesically Convex Functions}

\author{Zhanyuan Cai, Emre Sahinoglu and Shahin Shahrampour
\thanks{This work is supported in part by NSF Award ECCS-2240788   as well as NSF CAREER Award ECCS-2442321.}
\thanks{Z. Cai, E. Sahinoglu and S. Shahrampour are with the Department of Mechanical \& Industrial Engineering at Northeastern University, Boston, MA 02115, USA. 
        {\tt\small emails:\{cai.zha,sahinoglu.m, s.shahrampour\}@northeastern.edu}.}
}
\maketitle

\begin{abstract}
\everymath{}%
We study decentralized online optimization for strongly geodesically convex (strongly g-convex) losses on Riemannian manifolds with bounded sectional curvature, including positively curved manifolds. In centralized Riemannian optimization, strong g-convexity tightens the optimal regret from $\mathcal{O}(\sqrt{T})$ to $\mathcal{O}(\log T)$, where $T$ is the time horizon; in the decentralized Riemannian setting, however, existing methods address only g-convex losses, leaving the strongly g-convex regime unexplored. One challenge is that the required decaying step size in the centralized regime is incompatible with existing network-error analyses, which typically assume a fixed step size. First, we provide a general network-error analysis for time-varying schedules. Next, we build on this analysis to establish the first  $\mathcal{O}(\log T)$ static regret bound for decentralized online Riemannian gradient descent, matching the minimax-optimal rate for strongly-convex Euclidean online optimization. Finally, we prove the same  $\mathcal{O}(\log T)$ regret bound for the two-point bandit feedback setting using novel strong subconvexity arguments for the smoothed versions of the loss functions. 
\end{abstract}


\IEEEpeerreviewmaketitle

\section{Introduction}
Distributed online optimization has been studied extensively in the control literature with applications in multi-robot target tracking, distributed resource allocation and distributed economic dispatch. In this setting, a network of $n$ agents collaboratively minimizes a time-varying {\it global} loss $f_t = \tfrac{1}{n}\sum_{i=1}^{n} f_{i,t}$ while each agent $i$ observes only its {\it local} loss $f_{i,t}$ sequentially and can exchange information with neighboring agents. The performance of distributed online algorithms is measured by \emph{regret}, the difference between the collective cumulative loss of agents and that of the best single decision in hindsight. 

In Euclidean spaces, consensus-based gradient methods effectively solve this problem and achieve sublinear regret bounds~\cite{HosseiniChapmanMesbahi2016, ShahrampourJadbabaie2018, AkbariGharesifardLinder2017,YanSundaramVishwanathanQi2013}. These guarantees, however, rest on the vector-space structure of Euclidean space and do not extend readily to curved manifolds on which many modern engineering problems are naturally posed. Such manifold optimization problems are prevalent: collaborative geometric estimation for SLAM~\cite{tian2022distributed}, distributed PCA on orthogonality-constrained manifolds~\cite{gang2022linearly}, and consensus-based multi-robot active mapping over probability-simplex and $SE(3)$ variables~\cite{asgharivaskasi2025riemannian}.

In these manifold optimization problems, the Euclidean consensus update is inapplicable, since weighted averaging requires a vector-space structure; \emph{intrinsic} methods replace it with a tangent-space approximation of the weighted Fr\'echet mean~\cite{SarletteSepulchre2009, TronAfsariVidal2013,nguyen2026intrinsic,wang2025distributed}, whereas \emph{extrinsic} methods exploit an ambient Euclidean embedding and so apply only to a narrow class of manifolds, such as the Stiefel and other compact submanifolds~\cite{ChenGarciaHongShahrampour2021, ChenGarciaHongShahrampour2023, DengHu2025}. Both lines of work, however, target offline problems with non-convex or g-convex objectives.

In both Euclidean online optimization and centralized online Riemannian optimization, strong convexity sharpens the attainable regret from $\mathcal{O}(\sqrt{T})$ to $\mathcal{O}(\log T)$~\cite{HazanAgarwalKale2007, WangTuHongWuShi2023}, a rate that is minimax-optimal for strongly convex losses~\cite{abernethy2008optimal}. In the {\it decentralized} online Riemannian setting, however, this improvement remains out of reach: the only existing algorithms \cite{ChenSun2024,SahinogluShahrampour2025Decentralized} attain $\mathcal{O}(\sqrt{T})$ regret for g-convex losses and do not exploit strong g-convexity. This raises a basic question:
\begin{tcolorbox}[colback=gray!5, colframe=black!60, sharp corners, boxrule=0.5pt, left=6pt, right=6pt, top=4pt, bottom=4pt] \centering\emph{Can the logarithmic regret improvement from strong convexity be attained in decentralized online Riemannian optimization on
manifolds?}
\end{tcolorbox}
We answer the above question affirmatively by resolving three main challenges. First, network error: the ingredient delivering the logarithmic rate, a decaying step size, is incompatible with the existing analyses. All prior consensus analyses in decentralized Riemannian optimization assume a fixed step size~\cite{ChenSun2024,SahinogluShahrampour2025Decentralized}, under which the per-step error reduces to a stationary term. With a decaying step size, the iterates are re-perturbed by a different amount at every round, so the network error becomes a time-weighted accumulation of past step sizes, and the static argument no longer applies. Second, positive curvature: the metric projection onto the feasible set may not be non-expansive, which introduces an extra term in the optimization error; controlling it requires a shifted step size $\eta_t = O\!\big(1/(t+c_0)\big)$ rather than the standard $O(1/t)$~\cite{WangTuHongWuShi2023}. Third, the bandit feedback setting when the gradients are unavailable only compounds both effects: smoothing replaces each loss by a surrogate and degrades strong g-convexity, and keeping the queried points feasible shrinks the domain, so the smoothing and projection errors incorporate into the optimization and network terms simultaneously. Preserving the logarithmic rate therefore demands that a smoothing parameter be tuned jointly with the decaying schedule, so that each of these errors contributes only a constant. 

\textit{Contributions:} We study decentralized online Riemannian optimization for $\mu$-strongly g-convex losses on manifolds with bounded sectional curvature, including positively curved manifolds, and establish $\mathcal{O}(\log T)$ static regret under both full-information and two-point bandit feedback settings. Our work focuses on the curvature-aware consensus step and the two-point gradient estimator of~\cite{SahinogluShahrampour2025Decentralized}; the main contribution is the analysis of the strongly g-convex regime with a decaying step size to establish a novel logarithmic regret. Concretely:

\begin{itemize}
\item \textbf{First $\mathcal{O}(\log T)$ regret in decentralized online Riemannian optimization.} We establish the first $\mathcal{O}(\log T)$ static regret bound for strongly g-convex losses in the decentralized online Riemannian setting, on manifolds of bounded curvature including positively curved ones (Theorem \ref{thmregretboundfull}). The bound matches, up to curvature-dependent constants, the optimal logarithmic rate for strongly convex losses in both Euclidean optimization and centralized Riemannian optimization, and its $(1-\sigma_2(W))^{-1}$ dependence on the spectral gap of the mixing matrix $W$ is the expected price of decentralization.
\item \textbf{Consensus and optimization error analysis under a decaying step size.} Achieving this rate requires a decaying step size, which the existing fixed-step-size analyses \cite{ChenSun2024, SahinogluShahrampour2025Decentralized} do not cover. We develop a curvature-aware network error analysis valid for decaying schedules, in which the error accumulates as a time-weighted sum of past step sizes rather than a single stationary term, and we show that this sum remains $\mathcal{O}(\log T)$. On positively curved manifolds, where the metric projection may not be non-expansive, the same shifted decaying schedule also keeps the optimization error term controlled.
\item \textbf{Computationally efficient two-point bandit feedback.} We extend the analysis to the two-point bandit feedback, where each agent observes only two function values per round, using the estimator of~\cite{SahinogluShahrampour2025Decentralized,SahinogluSunShahrampour2025finite}. Since smoothing does not preserve strong g-convexity, we develop a novel strong sub-g-convexity argument that recovers a strong-convexity-type inequality for the smoothed surrogate up to a controllable slack (Lemma \ref{musubconvlem}), which yields the logarithmic rate (Theorem \ref{thmbanditreg}). The tighter $\mathcal{O}(\log T)$ target also forces a smaller smoothing radius than the g-convex case, $\delta = O(T^{-2})$ in place of $O(T^{-1})$, chosen jointly with the decaying schedule so that the smoothing, shrinkage, and projection errors remain constant and the $\mathcal{O}(\log T)$ rate is preserved.
\end{itemize}

\begin{table*}[t]

\caption{\centering Static regret bounds for online geodesically convex optimization on Riemannian manifolds.
$^\dagger$: separation oracle; $^\ddagger$: linear optimization oracle; $^{\dagger\dagger}$: static case;
cent: centralized; dec: decentralized; g-cvx: g-convex; sg-cvx: strongly g-convex}

\label{sampletable}
\begin{center}
\begin{small}
\begin{sc}
\begin{tabular}{C{4.9cm}ccccc}
\toprule
Reference & Objective & Manifold & Setting & Feedback & Regret Bound\\
\midrule

Wang et al.~\cite{WangTuHongWuShi2023}
& g-cvx
& Riemannian
& cent
& Gradient / 2-BAN
& $\mathcal{O}(\sqrt{T})$ \\

Wang et al.~\cite{WangTuHongWuShi2023}
& sg-cvx
& Riemannian
& cent
& Gradient / 2-BAN
& $\mathcal{O}(\log T)$ \\

Hu et al.~\cite{HuWangAbernethy2023}
& g-cvx
& Hadamard
& cent
& Gradient / 2-BAN$^\dagger$
& $\mathcal{O}(\sqrt{T})$ \\

Hu et al.~\cite{HuWangAbernethy2023}
& sg-cvx
& Hadamard
& cent
& Gradient$^\ddagger$
& $\mathcal{O}(T^{2/3}\log T)$ \\

Chen \& Sun~\cite{ChenSun2024}
& g-cvx
& Hadamard
& dec
& Gradient
& $\mathcal{O}(\sqrt{T})^{\dagger\dagger}$ \\

Sahinoglu \& Shahrampour~\cite{SahinogluShahrampour2025Decentralized}
& g-cvx
& Riemannian
& dec
& Gradient / 2-BAN
& $\mathcal{O}(\sqrt{T})$ \\

\textbf{Our work}
& sg-cvx
& Riemannian
& dec
& Gradient / 2-BAN
& $\mathcal{O}(\log T)$ \\

\bottomrule
\end{tabular}
\end{sc}
\end{small}
\end{center}
\vskip -0.2in
\end{table*}

\subsection{Related Work}

\textit{(i) Online optimization in Euclidean spaces.}
Online convex optimization (OCO) is well understood in the Euclidean space. The convex case admits $\mathcal{O}(\sqrt{T})$ regret~\cite{Zinkevich2003}, and strong convexity improves this rate to $\mathcal{O}(\log T)$~\cite{HazanAgarwalKale2007}. These guarantees, together with their dynamic-regret refinements~\cite{jadbabaie2015online,Mokhtari2016,ChangShahrampour2021}, form the standard backdrop for OCO. Moving from the single-agent to multi-agent scenario, consensus-based online gradient and mirror descent methods recover sublinear regret while incurring a cost for consensus in terms of the spectral gap of the communication matrix~\cite{HosseiniChapmanMesbahi2016, ShahrampourJadbabaie2018, AkbariGharesifardLinder2017, YanSundaramVishwanathanQi2013}. Recent analyses further attain nearly optimal rates for both convex and strongly convex losses~\cite{WanTuZhang2024}. The picture is essentially complete in the Euclidean setting, encompassing both the convex and strongly convex regimes as well as the centralized and decentralized settings.

\textit{(ii) Riemannian optimization.} Centralized optimization on manifolds has been well developed over the past two decades~\cite{AbsilMahonySepulchre2009, Sato2021}: global complexity guarantees for g-convex objectives, including the strongly g-convex case, are now available~\cite{ZhangSra2016}, together with accelerated~\cite{KimYang2022, MartinezRubio2022} and stochastic~\cite{Bonnabel2013} variants. The decentralized setting, however, has not kept pace and is broadly divided according to how consensus is enforced. \emph{Intrinsic} methods average through the tangent-space Fr\'echet-mean approximation and so apply on general manifolds with bounded curvature~\cite{ChenSun2024, wang2025distributed,nguyen2026intrinsic}, whereas \emph{extrinsic} methods rely on an ambient Euclidean embedding and are therefore confined to specific embedded manifolds, such as the Stiefel manifold~\cite{ChenGarciaHongShahrampour2021, ChenGarciaHongShahrampour2023} and compact submanifolds~\cite{DengHu2025, chen2025decentralized}. Nevertheless, prior analysis focuses on non-convex or g-convex objectives; {\it strong g-convexity}, the property that leads to the logarithmic regret in online optimization, has {\it not} been exploited in any decentralized Riemannian method.

\textit{(iii) Online Riemannian optimization.} For the centralized case, $\mathcal{O}(\sqrt{T})$ and $\mathcal{O}(\log T)$ regret rates are available for g-convex and strongly g-convex losses, respectively, with a matching lower bound in the g-convex case and an analysis that reaches beyond Hadamard manifolds to positively curved domains~\cite{WangTuHongWuShi2023}. A growing set of refinements has also followed, including zeroth-order~\cite{Maass2022}, optimistic~\cite{WangYuanHongHuWangShi2025}, projection-free~\cite{HuWangAbernethy2023}, and curvature-independent~\cite{SahinogluShahrampour2025Horospherical} methods, while the two-point bandit feedback we adopt traces to the multi-point and gradient-free estimators of Euclidean online optimization~\cite{AgarwalDekelXiao2010,FlaxmanKalaiMcMahan2005 }. For the {\it decentralized} online setting, however, only two prior works exist~\cite{ChenSun2024, SahinogluShahrampour2025Decentralized}, both restricted to g-convex losses and fixed step sizes. Strong g-convexity requires a decaying schedule, not addressed by any of these analyses. We close this gap by developing the first regret analysis for strongly g-convex losses in the decentralized regime.

\section{Preliminary}
\subsection{Background on Riemannian Optimization}
\label{subsec:prelim}

We consider a decentralized online optimization problem, defined over a $d$-dimensional complete Riemannian manifold $\mathcal{M}$, equipped with a Riemannian metric $g$. For any point $x \in \mathcal{M}$, the tangent space is denoted by $T_x\mathcal{M}$. $\mathcal{B}_{T_{x}\mathcal{M}}(r)$ (respectively,  $\mathcal{S}_{T_{x}\mathcal{M}}(r)$) denotes the ball (respectively, sphere) in tangent space $T_x\mathcal{M}$ centered at the origin with radius $r\ge0$.

A curve on $\mathcal{M}$ is called a \emph{geodesic} if it locally minimizes length, playing the role of a straight line in Euclidean space. The Riemannian metric $g$ induces a smoothly varying inner product\footnote{The subscript $x$ is omitted when clear from the context.} $\langle \cdot, \cdot \rangle_x$ and a geodesic distance function $d(x, y) \coloneqq \inf_{\gamma} \int_{0}^1\norm{ \gamma'(t)} dt$, where the infimum is over piecewise smooth curves $\gamma$ with $\gamma(0)=x$ and $\gamma(1)=y$. For a point $x \in \mathcal{M}$ with tangent space $T_x\mathcal{M}$, the exponential map $\Exp_x : T_x\mathcal{M} \to \mathcal{M}$ is defined as $\gamma(t)=\Exp_x(tu)$ for a vector $u \in T_x\mathcal{M}$, such that $\gamma(0) = x$ and $\gamma'(0) = u$.

The sectional curvature at a point $x$ of $\mathcal{M}$ measures how the manifold bends within a two-dimensional plane of directions through $x$. We say $\mathcal{M}$ has bounded sectional curvature if this is uniformly bounded across all such planes and all points of $\mathcal{M}$. Manifolds with everywhere non-positive sectional curvature are called Hadamard manifolds; in this paper, we place no restriction on the sign of the curvature and allow it to be positive. A key consequence of positive curvature is that $\Exp_x$ need not be a global diffeomorphism: it is invertible only within the injectivity radius $r_{\mathrm{inj}}>0$, the largest radius on which $\Exp_x$ remains a diffeomorphism.

A subset $\mathcal{X} \subseteq \mathcal{M}$ is geodesically convex (g-convex) if, any $x, y \in \mathcal{X}$ are joined by a geodesic $\gamma \subset \mathcal{X}$, and uniquely g-convex if this geodesic is unique. Since $\Exp_x$ is locally a diffeomorphism within $r_{\mathrm{inj}}$, it admits an inverse $\Log_x :\mathcal{M} \to T_x\mathcal{M}$ on a uniquely g-convex set.

For a differentiable function $f : \mathcal{M} \to \mathbb{R}$, the \textit{Riemannian gradient} $\grad f(x) \in T_x\mathcal{M}$ is the unique tangent vector satisfying
\begin{equation*}
\langle \grad f(x), v \rangle_x = Df(x)[v], \quad \forall v \in T_x\mathcal{M},    
\end{equation*}
where $Df(x)[v]$ denotes the directional derivative of $f$ at $x$ along $v$ \cite{boumal2023introduction}. A function $f : \mathcal{M} \to \mathbb{R}$ is said to be $\mu$-strongly g-convex ($\mu \geq 0$) if it is $\mu$-strongly convex along every geodesic; the case $\mu=0$ recovers g-convexity. For differentiable functions, this is equivalent to
\begin{equation*}
f(y)\ge f(x)+\langle\mathrm{grad}\,f(x),\mathrm{Log}_x(y)\rangle_x+\tfrac{\mu}{2}\,d^2(x,y),    
\end{equation*} 
$\forall x, y \in \mathcal{M}$. This generalizes Euclidean strong convexity by requiring a uniform quadratic lower bound along geodesics.

\subsection{Network Model, Consensus Update, and Regret}

We consider a network of $n$ agents collaboratively minimizing a sequence of global objectives over a uniquely g-convex compact subset $\mathcal{X} \subseteq \mathcal{M}$. Communication among the agents is modeled by a doubly stochastic weight matrix $W = [w_{ij}]$, where $w_{ij}\ge0$, and the entries in each row and each column sum to 1. In particular, $w_{ij} > 0$ if agent $i$ can receive information from agent $j$. Otherwise, $w_{ij} = 0$.

The Euclidean consensus update $x_i = \sum_{j=1}^n w_{ij} y_j$ is inapplicable on a manifold since this weighted average may be infeasible due to nonlinearity of the manifold. Its intrinsic analogue is the weighted Fr\'echet mean $x_i = \arg\min_{y \in \mathcal{X}} \{\sum_{j=1}^n w_{ij} d^2(y, y_j)\}$, which requires an auxiliary optimization at every step, so instead, we use the first-order approximation
\begin{equation*}
x_i=\mathrm{Exp}_{y_i}\Big(s\sum_{j=1}^n w_{ij}\,\mathrm{Log}_{y_i}(y_j)\Big),    
\end{equation*}
where the consensus step size $s>0$ is chosen according to the curvature of $\mathcal{M}$ \cite{ChenSun2024,SahinogluShahrampour2025Decentralized}. This moves $y_i$ along a tangent direction that depends on its neighbors.

At each round $t \in \{1, \dots, T\}$: (i) every agent $i$ selects $x_{i,t}\in\mathcal{X}$ using information from rounds $1,\dots,t-1$; (ii) the environment reveals a $\mu$-strongly g-convex loss $f_{i,t}: \mathcal{X} \to \mathbb{R}$; (iii) agent $i$ incurs a loss and receives some feedback. The feedback consists of the gradient $\mathrm{grad}\,f_{i,t}(x_{i,t})$ in the full information setting (Section \ref{sec:full_info}) and two function evaluations in the bandit setting (Section \ref{sec:bandit}). Note that agent $i$ has no knowledge of $f_{i,t}$ when choosing $x_{i,t}$.

The collective goal is to minimize the static regret relative to a fixed comparator $x^* \in \mathcal{X}$. Under full information feedback, regret is defined as
\begin{equation}\label{defregretfull}
    \Reg^{\mathrm{full}}(T)\coloneqq\frac{1}{n}\sum_{i=1}^n\sum_{t=1}^Tf_t(x_{i,t})-\sum_{t=1}^Tf_t(x^*),
\end{equation}
where $f_t(x) = \frac{1}{n}\sum_{i=1}^{n} f_{i,t}(x)$ is the global objective function at round $t$. This standard definition is commonly adopted in decentralized online optimization \cite{ChenSun2024,ShahrampourJadbabaie2018,SahinogluShahrampour2025Decentralized}.

The bandit regret definition differs from \eqref{defregretfull} only in the evaluation point. When the gradient is unavailable or costly to compute, we instead consider the two-point bandit setting, in which each agent $i$ observes the function values at the two queried points, $x_{i,t,1}$ and $x_{i,t,2}$, rather than $x_{i,t}$, and the regret is defined as 
\begin{equation}\label{defregretbandit}
\resizebox{1.03\columnwidth}{!}{$%
    \begin{aligned}
        \Reg^{\mathrm{bandit}}(T)\coloneqq&\frac{1}{2n}\sum_{i=1}^n\sum_{t=1}^Tf_t(x_{i,t,1})+f_t(x_{i,t,2})-\sum_{t=1}^Tf_t(x^*).
    \end{aligned}$}
\end{equation}
The two evaluations are used to build a two-point gradient estimator at $x_{i,t}$, as we will describe in Section \ref{sec:bandit}. This regret definition aligns with prior work on Riemannian online optimization \cite{WangTuHongWuShi2023} and its decentralized extension \cite{SahinogluShahrampour2025Decentralized}.

\subsection{Assumptions}

The theoretical analysis relies on the following standard assumptions on the network topology, the manifold geometry, and the properties of the loss functions.

\begin{assum}
\label{assum:network}
The network is connected and the communication matrix $W \in \mathbb{R}^{n \times n}$ is symmetric and doubly stochastic. $\sigma_2(W)$ denotes the second largest singular value of the matrix $W$, for which we have that $\sigma_2(W) \in [0,1)$.  
\end{assum}

This assumption is standard in the literature on decentralized online optimization \cite{ChenSun2024,SahinogluShahrampour2025Decentralized,ShahrampourJadbabaie2018}. Generally, a smaller $\sigma_2(W)$ corresponds to a better-connected network, facilitating faster information propagation among the agents.
\begin{assum}
    \label{assum:sectcurv}
    The sectional curvature $K$ on $\mathcal{M}$ is bounded, such that $K_{\min} \le K \le K_{\max}$. The diameter of set $\mathcal{X}$ is bounded by $D$. If $K_{\max}>0$, we further assume that $D < r_\mathrm{cx}$, where $r_\mathrm{cx}:=\frac{1}{2}\min\{r_\mathrm{inj},\frac{\pi}{\sqrt{K_{\max}}}\}$, and $r_\mathrm{inj}$ denotes the injectivity radius.  
\end{assum}

This assumption (1) implies that $\mathcal{X}$ is uniquely g-convex, so the Riemannian logarithm is well defined on $\mathcal{X}$, and (2) it rules out antipodal pairs, allowing meaningful strongly g-convex objectives on positively curved manifolds \cite{Yau1974Convex}. These geometric conditions are used in recent work on online Riemannian optimization beyond Hadamard manifolds \cite{WangTuHongWuShi2023,SahinogluShahrampour2025Decentralized}. 
\begin{assum}\label{assum:lossfunc}
    For all $i \in \{1, \dots, n\}$, the local losses $\{f_{i,t}\}_{t=1}^T$ are differentiable, $\mu$-strongly g-convex, and $L$-Lipschitz on the compact domain $\mathcal{X}$.
\end{assum}

Since $\mathcal{X}$ is compact and the losses are continuous, the sum $\sum_{t=1}^T f_t$ attains its minimum on $\mathcal{X}$. We fix the minimizer $x^*=\arg\min_{x\in\mathcal{X}}\sum_{t=1}^T f_t(x)$ as the comparator against which the static regrets \eqref{defregretfull} and \eqref{defregretbandit} are measured. The loss regularity in Assumption \ref{assum:lossfunc} is likewise standard in geodesically convex optimization \cite{ZhangSra2016,WangTuHongWuShi2023,SahinogluShahrampour2025Decentralized}. However, the strong g-convexity is what we specifically analyze in this paper.

\section{Decentralized Online Riemannian Optimization: Full Information Feedback}\label{sec:full_info}

In this section, we establish the first $\mathcal{O}(\log T)$ static regret bound for decentralized online Riemannian optimization in the full information setting (Theorem \ref{thmregretboundfull}). 

We analyze the same two-step algorithm as in \cite{SahinogluShahrampour2025Decentralized}, but under a decaying step size, as required to achieve the strongly g-convex regret rate. This change necessitates the new consensus analysis developed in Lemma  \ref{lemnetworkerror}. At each time step $t$, every agent $i$ performs the following two updates:
\begin{align}   y_{i,t+1}&=
\mathcal{P}_{\mathcal{X}}\!\left(\Exp_{x_{i,t}}\!\left(-\eta_t \grad_{i,t}\right)\right),\tag{I}\label{step1}\\ x_{i,t+1}&=\Exp_{y_{i,t+1}}\!\Big(s\sum_{j=1}^n w_{ij}\Log_{y_{i,t+1}}(y_{j,t+1})\Big),\tag{II}\label{step2}
\end{align}
where we write $\grad_{i,t} := \grad f_{i,t}(x_{i,t})$ for brevity, $\mathcal{P}_{\mathcal{X}}:\mathcal{M}\rightarrow\mathcal{X}$ denotes the  projection onto $\mathcal{X}$, defined by $\mathcal{P}_{\mathcal{X}}(x)\coloneqq \arg\min_{y\in\mathcal{X}}d(x,y)$, and $s>0$ is the consensus step size.

In our analysis, the regret decomposes into two components (errors) corresponding to the two steps of the update. Step~\eqref{step1} is a projected Riemannian gradient descent update: agent $i$ moves along the negative gradient via the exponential map and projects back onto $\mathcal{X}$, contributing to the \emph{online optimization error}. Step~\eqref{step2} is a curvature-aware consensus step: agent $i$ takes a step of size $s$ toward the weighted Fr\'echet mean of its neighbors' post-gradient iterates $\{y_{j,t+1}\}$ for any $j$ in the neighborhood of $i$, contributing to the \emph{network error}, i.e., the discrepancy among agents' iterates due to imperfect communication. Let us now discuss the {\it challenges} in the analysis of each component.

\textit{(a) Online optimization error:} The principal difficulty in this term is the metric projection $\mathcal{P}_{\mathcal{X}}$, whose non-expansiveness is not guaranteed on positively curved manifolds. On Hadamard manifolds, $\mathcal{P}_{\mathcal{X}}$ is non-expansive and the analysis proceeds as in the Euclidean case. On manifolds with positive sectional curvature, non-expansiveness may fail to hold, and the projection introduces an additional error term. To keep this term summable while preserving the $\mathcal{O}(\log T)$  rate, we adopt the shifted decreasing step size $\eta_t = O(1/(t+c_0))$ of \cite{WangTuHongWuShi2023}, where $c_0>0$ certifies $\eta_t L \leq D$ for all $t\geq 1$. 

\textit{(b) Network error:} The network error quantifies the deviation of each agent's iterate from the network Fr\'echet mean, and insufficient consensus degrades the regret. 
Our analysis also reveals an extra interference term that does not contribute to regret, provided $s$ is chosen appropriately with respect to the curvature (see \eqref{eq:interference} in Section \ref{app:thmfull}). Consequently, only the network error itself, namely the failure to reach exact consensus, contributes to the final regret bound. Unlike the online-optimization term, this is where decentralization and the decaying step size interact, and where the fixed-step-size analyses of \cite{ChenSun2024, SahinogluShahrampour2025Decentralized} do not apply.

We first bound the network error. The consensus step (\ref{step2}) contracts the Fr\'echet variance of the agents' iterates at a linear rate. Let us denote by $\bar{y}_t$ the Fr\'echet mean of $\{y_{i,t}\}_{i=1}^n$ and by $\bar{x}_t$ the Fr\'echet mean of $\{x_{i,t}\}_{i=1}^n$. \cite[Theorem III.2]{SahinogluShahrampour2025Decentralized} proves that $\sum_{i=1}^n d^2(x_{i,t},\bar{x}_t)\le\rho \sum_{i=1}^n d^2(y_{i,t},\bar{y}_t)$, where the factor $\rho\in(0,1)$ depends on the consensus step size $s$, the spectral gap $1-\sigma_2(W)$, and the curvature. This contraction is a property of the consensus step alone and is unaffected by the choice of gradient step size $\eta_t$. What the decaying $\eta_t$ does change is the disturbance the contraction must absorb: Step (\ref{step1}) re-perturbs the iterates by an $\eta_t$-dependent amount, so the network error behaves like a geometric convolution $\sum_k \rho^{(t-k)/2}\eta_k$ of past step sizes. Prior analyses \cite{ChenSun2024,SahinogluShahrampour2025Decentralized} evaluate this sum only for a constant $\eta
$, for which it collapses to a single static term. The decaying schedule requires summing the convolution against a time-varying sequence. The following lemma carries out this analysis and shows the resulting sum remains $\mathcal{O}(\log T)$. 

\begin{lem}[Network error]\label{lemnetworkerror}
    Let Assumptions \ref{assum:network}, \ref{assum:sectcurv}, and \ref{assum:lossfunc} hold, and suppose all agents share a common initialization $x_{i,1}=x_1 \in \mathcal{X}$. If we set $s=(2C_1)^{-1}C_2$, then for any step size sequence $\{\eta_t \}$ with $\eta_t\le \tfrac{D}{L}$ and for all $t\ge2$, running Steps \eqref{step1} and \eqref{step2} gives
    \begin{equation}\label{networkerror}
        \begin{aligned}
            \frac{1}{n}\sum_{i=1}^nd(x_{i,t},\bar{x}_t)\le L\cdot\sum_{k=1}^{t-1}\left(\rho^{\frac{t-k}{2}} \eta_k\right),
        \end{aligned}
    \end{equation}
    where $\rho := 1 - \frac{C_2^2 (1 - \sigma_2(W))}{2 C_1 (1 + C_3 D^2)^2}$, and $C_1, C_2,C_3$ are curvature-dependent constants defined in Appendix \ref{subsecconst}. In particular, under $\eta_t=(\mu t+L/D)^{-1}$, for all $T\ge1$, it holds that
    \begin{equation*}
        \frac{1}{n}\sum_{t=1}^T\sum_{i=1}^nd(x_{i,t},\bar{x}_t)\le \frac{2\,L}{\mu(1-\rho)}\log\left(1+\frac{T\mu D}{L}\right).
    \end{equation*}
\end{lem}

Combining the network-error bound of Lemma \ref{lemnetworkerror} with the online optimization error analysis, we obtain the following static regret bound.

\begin{thm}\label{thmregretboundfull}
    Let Assumptions \ref{assum:network}, \ref{assum:sectcurv}, and \ref{assum:lossfunc} hold, and suppose all agents share a common initialization $x_{i,1}=x_1 \in \mathcal{X}$.  Running the two Steps \eqref{step1} and \eqref{step2} for $T$ iterations with step size $\eta_t=(\mu t+L/D)^{-1}$ and $s=(2C_1)^{-1}C_2$, we have the following static regret bound:
    \begin{equation*}
    \begin{aligned}
        \Reg^{\mathrm{full}}(T)
        &\le R_1\log \left(1+\frac{T\mu D}{L}\right)+R_2,
    \end{aligned}
    \end{equation*}
    where $R_1$ and $R_2$ are constants independent of $T$. See Appendix \ref{subsecconst} for their values.
\end{thm}

\begin{rem}
In the Euclidean decentralized online optimization setting, the static regret bound for strongly convex functions under the standard communication scheme is of order $\log T/(1-\sigma_2(W))$~\cite{YanSundaramVishwanathanQi2013}. Our bound matches this order in both $T$ and $\sigma_2(W)$ but has extra geometry-dependent constants reflecting the cost of optimization on a non-flat manifold.
\end{rem}

\begin{rem}
In the centralized online Riemannian optimization setting, an $\mathcal{O}(\log T)$ regret bound is established for strongly g-convex functions \cite{WangTuHongWuShi2023}. Our result recovers the same rate in the decentralized setting, with the additional cost of $(1-\sigma_2(W))^{-1}$ for consensus in multi-agent networks.
\end{rem}

\section{Decentralized Online Riemannian Optimization: Two-Point Bandit Feedback}
\label{sec:bandit}

Under two-point bandit feedback, agents observe only function values and must rely on zeroth-order information. Each agent queries its loss at two nearby points and forms a two-point gradient estimator in place of the Riemannian gradient. Establishing an $\mathcal{O}(\log T)$ bound here raises three difficulties beyond the full-information case. First, the estimator is associated with a smoothed surrogate of the loss, and smoothing does not preserve strong g-convexity. We should quantify the residual as a strong-sub-g-convexity property (Lemma \ref{musubconvlem}). Second, constraint handling introduces shrinkage and projection errors. Third, these terms must be controlled tightly enough to keep the logarithmic rate. With an appropriate smoothing parameter each contributes only a constant, so the $\mathcal{O}(\log T)$ bound carries over.

Let $x_{i,t}$ denote the decision of agent $i$ at time $t$. Following
\cite{SahinogluShahrampour2025Decentralized}, we sample a direction $u_{i,t}$ uniformly from the unit sphere
$\mathcal{S}_{T_{x_{i,t}}\mathcal{M}}(1)$ and form  two perturbed points
$x_{i,t,1}\coloneqq\Exp_{x_{i,t}}( \delta u_{i,t})$ and $x_{i,t,2}\coloneqq\Exp_{x_{i,t}}(-\delta  u_{i,t})$, at most $2\delta$
apart along a common geodesic through $x_{i,t}$. From their central
difference we define the estimator 
\begin{equation}\label{eq:gradestimator}
    \esgrad_{i,t}=\frac{d}{2\delta}\left(f_{i,t}(x_{i,t,1})-f_{i,t}(x_{i,t,2})\right) u_{i,t},
\end{equation}
where $\delta > 0$ is the smoothing parameter and $d$ is the dimension of $\mathcal{M}$.

First, we address the feasibility of the queries. Because the two-point estimator evaluates the loss at perturbed points $x_{i,t,1}$ and $x_{i,t,2}$, these points must remain in the feasible set $\mathcal{X}$. So, we constrain the decision variables $x_{i,t}$ to a shrinking subset of $\mathcal{X}$, which requires the following geometric assumption.

\begin{assum}\label{assum:domain-sandwich}
    There exists an interior reference point $p \in \mathcal{X}$ and radii $0 < r \leq R$ such that $\mathcal{B}_p(r) \subseteq \mathcal{X} \subseteq \mathcal{B}_p(R)$, where $\mathcal{B}_p(\rho)$ denotes the geodesic ball of radius $\rho$ centered at $p$.
\end{assum}

Following \cite{SahinogluShahrampour2025Decentralized}, restricting the iterates to the shrinking set $(1-\tau)\mathcal{X} \coloneqq \{\Exp_p((1-\tau)\Log_p(x)) \mid x \in \mathcal{X}\}$ with shrinkage factor $\tau \coloneqq \frac{\delta\theta}{r}$ guarantees that the perturbed points $x_{i,t,1}$ and $x_{i,t,2}$ remain in $\mathcal{X}$. The geometric constant $\theta>0$, which depends on the sectional-curvature bounds of $\mathcal{M}$ and the radii $R$ and $r$, is defined in \cite{SahinogluShahrampour2025Decentralized}.

Next, to analyze the estimator, we introduce the {\it smoothed} surrogate of each local loss as follows, 
\begin{equation*}
    \hat{f}^\delta_{i,t}(x)\coloneqq \int_{\mathcal{B}_{T_{x}\mathcal{M}}(1)} f_{i,t}\left(\Exp_x(\delta u)\right)\ dp(u),
\end{equation*}
where $p$ denotes the uniform measure on $\mathcal{B}_{T_{x}\mathcal{M}}(1)$. Since the update is driven by $\esgrad_{i,t}$, the algorithm effectively optimizes these surrogates rather than the losses themselves, so we analyze the regret against $\hat{f}^{\delta}_{i,t}$ and control the gap to $f_{i,t}$ via Lipschitzness. The crux is that $\hat{f}^{\delta}_{i,t}$ satisfies a strong-sub-g-convexity inequality with respect to $\esgrad_{i,t}$ up to a controllable slack (Lemma \ref{musubconvlem}), replacing the exact
strong g-convexity of the full-information case. We begin by formalizing the relaxed strong-convexity notion used in the analysis.
\begin{defn}
\label{def:subconvex}
A function $f:\mathcal{X} \subseteq \mathcal{M} \rightarrow \mathbb{R}$, together with a vector field $G$ on $\mathcal{X}$, is said to be $\mu$-strong $\lambda$-sub-g-convex for constants $\mu$, $\lambda\geq 0$ if, for all $x,y \in \mathcal{X}$, 
\begin{equation*}
    f(y)-f(x)-\langle G(x) ,\Log_x(y) \rangle - \tfrac{\mu}{2}d^2(x,y) \geq -\lambda.
\end{equation*}
\end{defn}

When $G = \grad f$, this recovers $\mu$-strong g-convexity up to the additive slack $\lambda$. The  above definition allows $G$ to be a surrogate for the gradient, such as the estimator $\esgrad$.

The next lemma shows that the expectation of the estimator $\esgrad_{i,t}$ satisfies this inequality for the smoothed surrogate $\hat{f}^{\delta}_{i,t}$, with a slack proportional to $\delta$.

\begin{lem}\label{musubconvlem}
    Suppose Assumption \ref{assum:sectcurv} holds, and let $f:\mathcal{X}\subseteq \mathcal{M} \rightarrow \mathbb{R}$ be $\mu$-strongly g-convex and $L$-Lipschitz on the uniquely g-convex domain $\mathcal{X}$. Fix a smoothing radius $0<\delta \le D$, let $\hat{f}^{\delta}(x)=\int_{\mathcal{B}_{T_{x}\mathcal{M}}(1)} f(\Exp_x(\delta u)) dp(u)$ be the smoothed surrogate, and define the gradient estimator
    \begin{equation}
    \label{eq:gdelta}
        \esgrad (x) = \frac{d}{2\delta} (f(\Exp_x(\delta u)) - f(\Exp_x(-\delta u)))u,
    \end{equation}
    where $u$ is uniformly distributed on $\mathcal{S}_{T_{x}\mathcal{M}}(1)$. Then, the pair $(\hat{f}^{\delta},G)$ with certifying field $G(x)\coloneqq \mathbb{E}_u [\esgrad(x)] \in T_{x}\mathcal{M}$ is $\mu$-strong $\delta( L C_4+2\mu D)$-sub-g-convex in the sense of Definition \ref{def:subconvex} for all $x,y\in (1-\tau)\mathcal{X}$, i.e.,  
    \begin{equation}\label{subconvex}
        \begin{aligned}
            \hat{f}^{\delta}(y)-\hat{f}^{\delta}(x) -\langle \mathbb{E}_u[\esgrad(x)],\Log_x(y) \rangle - \frac{\mu}{2}d^2(x,y)\\ 
            \geq -\delta ( L C_4+2\mu D),
        \end{aligned}
    \end{equation}
    where $C_4>0$ depends on $\kmin$, $\kmax$, and $D$ (see Appendix \ref{subsecconst}).
\end{lem}

\begin{rem}
The g-convex counterpart in \cite{SahinogluShahrampour2025Decentralized} establishes that $\hat{f}_{i,t}^{\delta}$ is $\delta L C_4$-sub-g-convex. Exploiting strong g-convexity, our lemma adds the quadratic term $\tfrac{\mu}{2}d^2(x,y)$ at the price of the extra $2\mu D$ slack. The same trade-off, a stronger convexity inequality for a larger additive constant, appears in centralized online Riemannian optimization \cite{WangTuHongWuShi2023}. And, when $\mu=0$, we recover the result of \cite{SahinogluShahrampour2025Decentralized} as expected.
\end{rem}

The update for each agent $i$ at time $t$ is then
\begin{align}
    y_{i,t+1} =&\ \mathcal{P}_{(1-\tau)\mathcal{X}}\left(\mathcal{P}_{\mathcal{X}}\left(\Exp_{x_{i,t}}\left(-\eta_{t}\esgrad_{i,t}\right)\right)\right),\tag{I'}\label{bandstep1'}
\end{align}
followed by the same consensus step \eqref{step2} as in the full-information setting. The bandit update thus differs from \eqref{step1}-\eqref{step2} only in two respects: the Riemannian gradient is replaced by the two-point estimator $\esgrad_{i,t}$ defined in \eqref{eq:gradestimator}, and the gradient iterate is additionally projected onto the shrinking set $(1-\tau)\mathcal{X}$.

We are now ready to establish the main theoretical guarantee regarding the regret bound in the two-point bandit setting.

\begin{thm}\label{thmbanditreg}
    Suppose that Assumptions \ref{assum:network},\ref{assum:sectcurv},\ref{assum:lossfunc} as well as \ref{assum:domain-sandwich} hold. Let the gradient estimator $\esgrad_{i,t}$ be generated as in \eqref{eq:gradestimator}. With step size $\eta_t = \left(\mu t + \frac{dL}{D}\right)^{-1}$, the smoothing parameter $\delta = T^{-2}$, and the shrinkage factor $\tau = \delta\theta/r$, running the update steps (\ref{bandstep1'}) and (\ref{step2}) for $T$ rounds yields the following expected static regret bound:
    \begin{equation*}
        \begin{aligned}
            \mathbb{E}[\Reg^{\mathrm{bandit}}(T)]
        \le& R_1'\log\!\Big(1+\frac{T\mu D}{dL}\Big) + R_2',
        \end{aligned}
    \end{equation*}
    where $R_1'$ and $R_2'$ are explicitly defined constants that are independent of the time horizon $T$ (see Appendix \ref{subsecconst}).
\end{thm}

Thus, although smoothing degrades strong g-convexity and introduces a gap between $f_{i,t}$ and $\hat{f}^{\delta}_{i,t}$, and the sampling constraint shrinks the feasible set, the $\mathcal{O}(\log T)$ rate is preserved for $\mu$-strongly g-convex losses on manifolds. In the centralized case ($n=1$), the theorem recovers the rate of \cite{WangTuHongWuShi2023}, and the extra $1-\sigma_2(W)$ dependence reflects the cost of decentralization.

\section{Numerical Experiments}\label{sec:experiments}

In this section, we conduct experiments to evaluate the performance of our algorithms. The first experiment is on a positively curved manifold, $\mathbb{S}^{7}\subset\mathbb{R}^{8}$, which poses greater geometric challenges than Hadamard manifolds due to possible projection errors. The second experiment is with a real-world dataset on the symmetric positive-definite (SPD) matrix manifold, which lets us analyze the decaying step size for strongly g-convex objectives on a Hadamard manifold 
\footnote{The code is available at: \href{https://github.com/ZhanyCai-opt/decentralized-online-riemannian-strongly-g-convex}{https://github.com/ZhanyCai-opt/decentralized-online-riemannian-strongly-g-convex}.}.

\subsection{Hyper-sphere experiment}

Adapting the online Fr\'echet-mean setup of \cite{SahinogluShahrampour2025Decentralized}, each agent $i$ observes the loss $f_{i,t}(x) = d^2(x, z_{i,t})$, with global objective $f_t(x) = \frac{1}{n}\sum_{i=1}^{n} f_{i,t}(x)$. The feasible set $\mathcal{X}$ is a geodesic ball of radius $7\pi/32$, so its diameter satisfies $D = 7\pi/16 < \pi/(2\sqrt{K_{\max}}) = \pi/2$ and Assumption \ref{assum:sectcurv} holds. By Hessian comparison \cite{petersen2006riemannian}, $d^2(\cdot, z)$ is $\mu$-strongly g-convex on $\mathcal{X}$.

\begin{figure}[t]
    \centering
    \includegraphics[width=\linewidth]{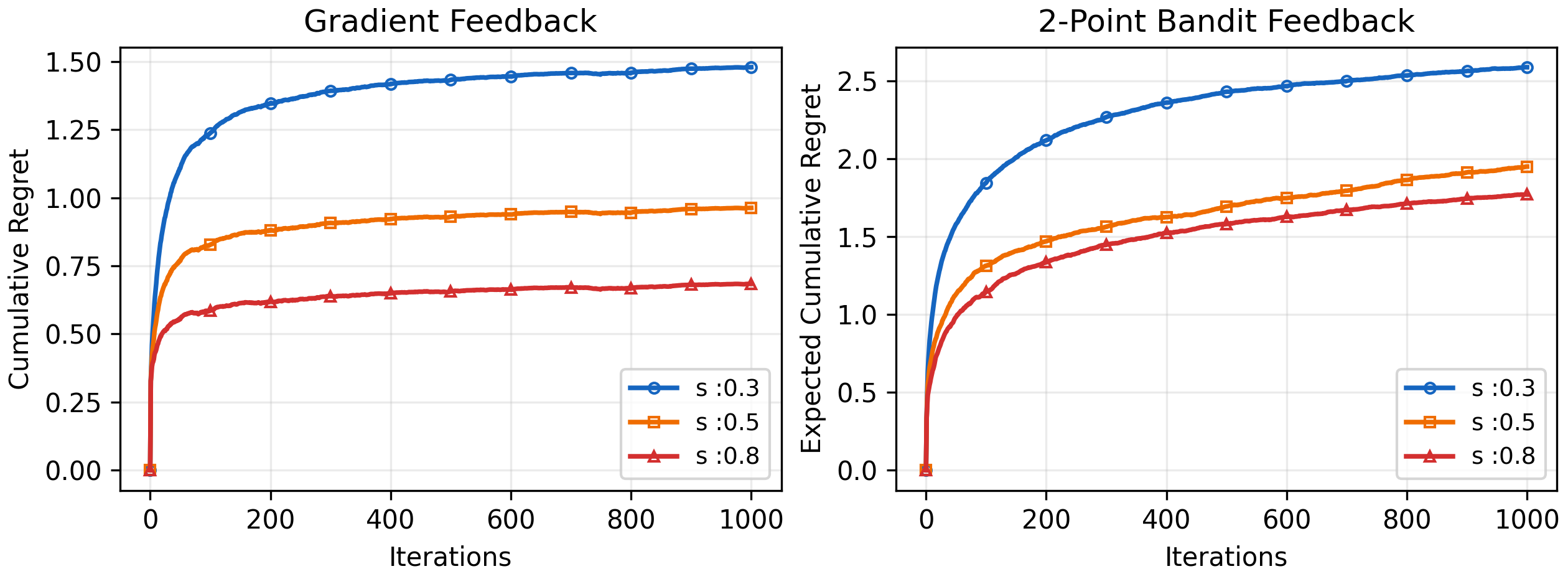}
    \caption{Static regret for strongly g-convex losses with consensus step size $s\in\{0.3,0.5,0.8\}$.}
    \label{fig:regret_s}
\end{figure}

\begin{figure}[t]
    \centering
    \includegraphics[width=\linewidth]{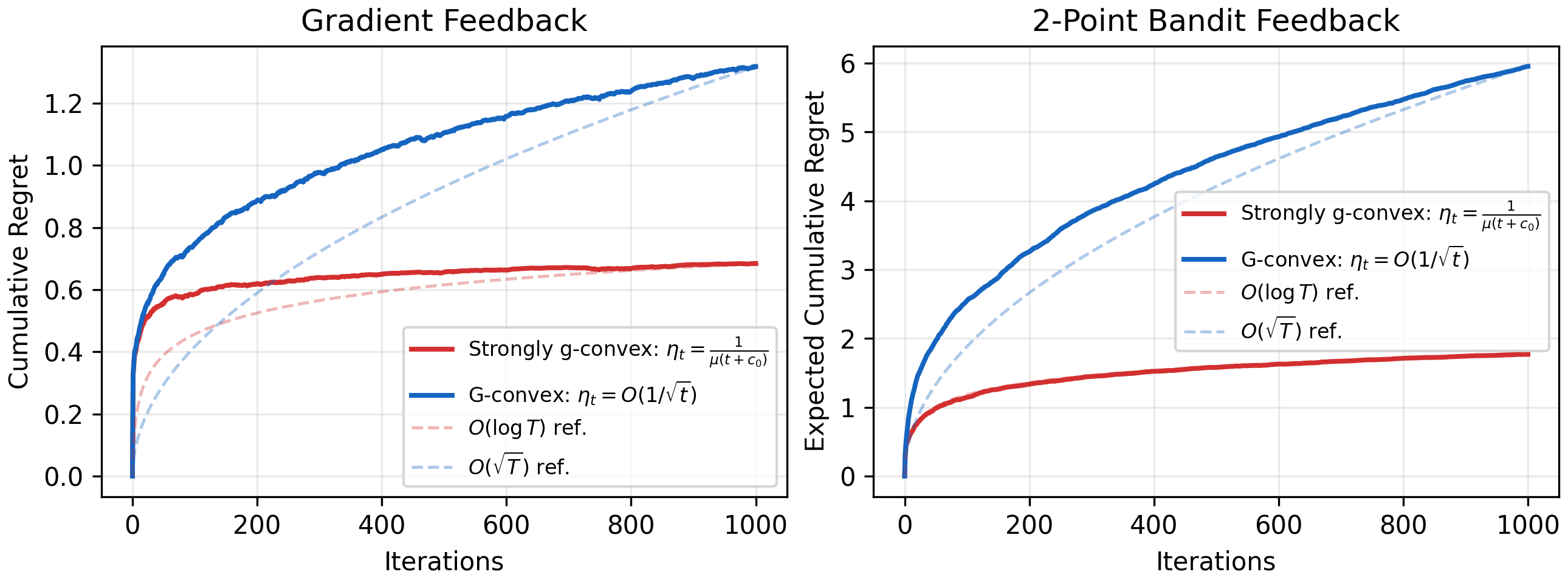}
    \caption{Strongly g-convex schedule, $O(\log T)$, versus g-convex baseline, $O(\sqrt{T})$, with $s=0.8$.}
    \label{fig:sc_vs_gc}
\end{figure}

We use $n = 10$ agents on a ring graph, where each agent communicates with its $k = 2$ nearest neighbors on each side, giving $\sigma_2(W) \approx 0.65$. The base points $\{z_i\}_{i=1}^{n}$ are sampled from an inner ball, and each $z_{i,t}$ is drawn uniformly from a small ball around $z_i$. We run both feedback models for $T = 1000$ rounds, using smoothing $\delta = T^{-2}$ under bandit feedback, and average the bandit results over $5$ Monte-Carlo simulation runs. Fig.~\ref{fig:regret_s} reports the static regret for consensus step-sizes $s \in \{0.3, 0.5, 0.8\}$ under full-information and two-point bandit feedback settings. Fig.~\ref{fig:sc_vs_gc} compares the strongly g-convex schedule against the g-convex baseline \cite{ChenSun2024,SahinogluShahrampour2025Decentralized} (with $\eta_t = O(1/\sqrt{t})$ and $\delta = O(1/T)$), together with $O(\log T)$ and $O(\sqrt{T})$ reference curves.

Fig.~\ref{fig:regret_s} shows that a larger consensus step size $s$ consistently reduces regret, and that bandit feedback incurs larger regret than full-information feedback, reflecting the estimator's variance and dimension dependence. Fig.~\ref{fig:sc_vs_gc} confirms the predicted separation: the strongly g-convex schedule grows logarithmically and quickly stabilizes, while the g-convex baseline keeps growing at the $O(\sqrt{T})$ rate under both feedback models.

\subsection{FLUXNET2015 experiment on SPD manifold}

We evaluate the proposed method on a real dataset, namely FLUXNET2015. We treat each flux-tower site as an agent and, for each week from 2010 to 2013, form an SPD matrix from the weekly correlation matrix of five meteorological variables. 
These matrices lie on the manifold $\mathrm{SPD(5)}$ equipped with the affine-invariant metric. Each agent's local loss is $f_{i,t}(X) = \tfrac{1}{2}\, d^2_{\mathrm{SPD}}(X, Z_{i,t})$, where $Z_{i,t}$ is the weekly SPD matrix at site $i$, and the iterates are constrained to a geodesic ball centered at the identity.

\begin{figure}[t]
    \centering
    \includegraphics[width=\linewidth]{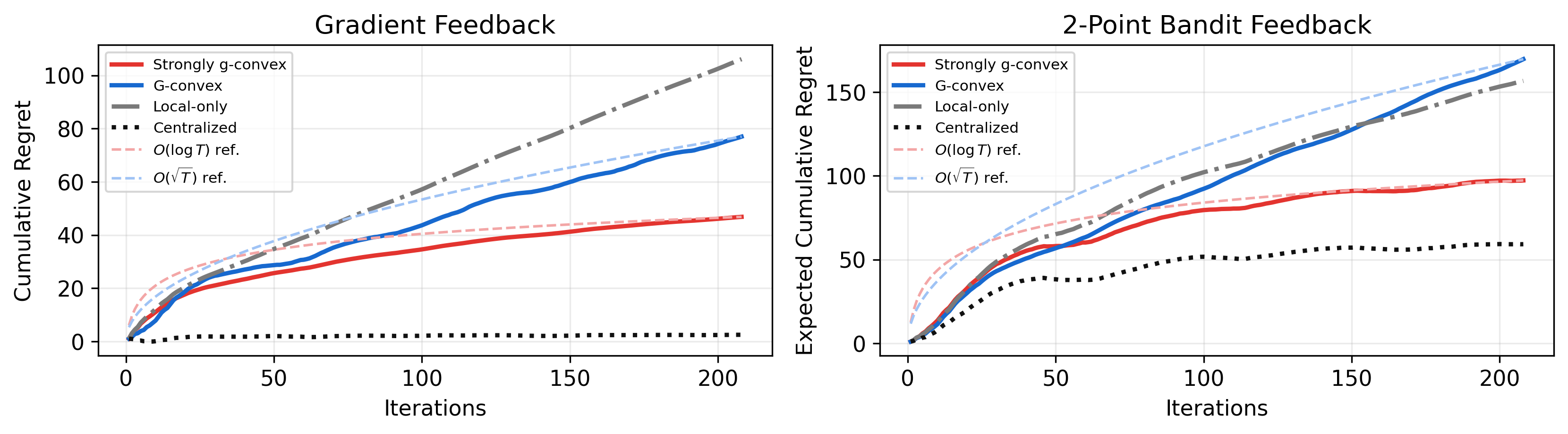}
    \caption{FLUXNET2015 results. Strongly g-convex regret, $O(\log T)$, versus g-convex baseline, $O(\sqrt{T})$.}
    \label{fig:fluxnet}
\end{figure}

We use $n = 20$ agents on a ring graph with three neighbors on each side and run the algorithms for $T = 208$ weeks. Since $\mathrm{SPD(5)}$ is Hadamard, the projection is non-expansive and the shift $c_0 = 0$ is sufficient. We compare the strongly g-convex schedule against a convex-stepsize baseline, a local-only baseline, the centralized method, and the two-point bandit variant. Fig.~\ref{fig:fluxnet} depicts the cumulative static regret, and the network consensus error with final regret values are tabulated in Table \ref{tab:fluxnet}.

Among the decentralized full-information methods, the strongly g-convex schedule achieves the lowest regret, $46.83$ versus $76.97$ for the convex baseline, confirming the advantage of the decaying stepsize on a Hadamard manifold. The local-only baseline is the worst, highlighting the necessity of communication, while the two-point bandit method incurs larger regret, as expected from zeroth-order estimation, but still grows sublinearly. The strongly g-convex method also attains a smaller consensus error than the convex baseline, indicating more agreement
among agents.
\begin{table}[t]
    \centering
    \caption{Final FLUXNET2015 results.}
    \label{tab:fluxnet}
    \begin{tabular}{l c c}
        \toprule
        \textsc{Algorithm} & \textsc{Final regret} & \textsc{Consensus error} \\
        \midrule
        \textsc{Strongly g-convex full} & $46.828497$ & $1.486785 \times 10^{-1}$ \\
        \textsc{G-convex full}            & $76.971451$ & $7.190531 \times 10^{-1}$ \\
        \textsc{Local-only full}        &$106.078126$ & $8.304925 \times 10^{-1}$ \\
        \textsc{Centralized full} & $2.522229$ & $0.000000$ \\ \hline
        \textsc{Strongly g-convex bandit}  & $97.297616$ & $1.602836 \times 10^{-1}$ \\
        \textsc{G-convex bandit}            & $169.713605$ & $1.682539$ \\
        \textsc{Local-only bandit}        & $156.744117$ & $8.835844 \times 10^{-1}$ \\ 
        \textsc{Centralized bandit} & $59.226032$ & $0.000000$ \\ 
        \bottomrule
    \end{tabular}
\end{table}

\section{Conclusion, Limitations, and Future Work}

We studied decentralized online Riemannian optimization for strongly g-convex losses on manifolds with bounded sectional curvature. Our analysis gives a network-consensus-error bound that holds for general time-varying step-size schedules under a mild condition (Lemma~\ref{lemnetworkerror}), and uses it to establish an $\mathcal{O}(\log T)$ regret bound under both full-information and two-point bandit feedback settings. This bound matches the lower bound for strongly convex losses in Euclidean online optimization, and is corroborated by our numerical experiments.

Two directions remain open. First, our bound carries the unaccelerated $(1-\sigma_2(W))^{-1}$ spectral-gap dependence of standard gossip; extending Euclidean accelerated-gossip schemes \cite{WanTuZhang2024, Wan2025BlackBox} to manifolds would potentially improve it. Second, dynamic regret under strong g-convexity is largely unexplored, and the nontrivial interaction of time-varying comparators with curvature makes this an interesting future direction.

\section{Appendix}
\subsection{Proof of Lemma \ref{lemnetworkerror}}\label{proofnetworkerror} 
By \cite[Theorem III.2]{SahinogluShahrampour2025Decentralized}, one application of the consensus step (\ref{step2}) contracts the Fr\'echet variance of the iterates as 
\begin{equation*}
    \sum_{i=1}^{n} d^{2}(x_{i, t+1}, \bar{x}_{t+1}) \leq \rho \sum_{i=1}^{n} d^{2}(y_{i, t+1}, \bar{y}_{t+1}), 
\end{equation*}
where the factor $\rho\in(0,1)$. Setting $a_t := \big(\sum_{i=1}^n d^2(x_{i,t},\bar x_t)\big)^{1/2}$, we then have 
\begin{align*}
a_{t+1} & \le  \sqrt{ \rho\sum_{i=1}^{n} d^{2}(y_{i, t+1}, \bar{y}_{t+1})}\\
&\leq \sqrt{\rho \sum_{i=1}^{n} d^{2}(y_{i, t+1}, \bar{x}_{t})} \\
&\leq \sqrt{\rho \sum_{i=1}^{n} d^{2}(x_{i, t}, \bar{x}_{t})} + \sqrt{\rho \sum_{i=1}^{n} d^{2}(x_{i, t}, y_{i, t+1})} \\
&\leq \sqrt{\rho \sum_{i=1}^{n} d^{2}(x_{i, t}, \bar{x}_{t})} + \sqrt{\rho n} \eta_t L. 
\end{align*}
where the second inequality uses the fact that the Fr\'echet mean $\bar{y}_{t+1}$ is the minimizer of the function $g(y)=\sum_i d^2(y,y_{i,t+1})$ over $\mathcal{X}$. The third follows from Minkowski's inequality applied to $d(y_{i,t+1},\bar x_t)\le d(x_{i,t},\bar x_t)+d(x_{i,t},y_{i,t+1})$, and the last uses the bound $d(x_{i,t},y_{i,t+1})\le \eta_t L$ due to Step~\eqref{step1}.

Under the common initialization $x_{i,1}=x_1$ we have $a_1=0$. Unrolling the recursion $a_{t+1}\le \sqrt{\rho}\,a_t + \sqrt{\rho n}\,\eta_t L$ then gives, for all $t\ge 2$, $a_t \le  \sqrt{n}\,L\sum_{k=1}^{t-1}\rho^{\frac{t-k}{2}}\eta_k$ and
\begin{align*}
\frac{a_t^2}{n} \;\le\; L^2\Big(\sum_{k=1}^{t-1}\rho^{\frac{t-k}{2}}\eta_k\Big)^2.
\end{align*}
By the Cauchy–Schwarz inequality, we have $\frac1n\sum_{i} d(x_{i,t},\bar{x}_t)\le \frac{1}{\sqrt{n}}\sqrt{\sum_i d^2(x_{i,t},\bar{x}_t)}$; hence,
\begin{equation*}
    \frac1n\sum_{i=1}^n d(x_{i,t},\bar x_t) \;\le\; L\sum_{k=1}^{t-1}\rho^{\frac{t-k}{2}}\eta_k.
\end{equation*}
Summing over $t$  and exchanging the order of summation,
\begin{equation*}
    \begin{aligned}
        \frac{1}{n}\sum_{t=1}^T\sum_{i=1}^nd(x_{i,t},\bar{x}_t)&\le L\cdot\sum_{t=2}^T\sum_{k=1}^{t-1}\left(\rho^{\frac{t-k}{2}} \eta_k\right)\\
        &\le L\cdot \sum_{k=1}^{T-1}\eta_k\sum_{t=k+1}^T\rho^{\frac{t-k}{2}}\\
        &\le \frac{L\cdot \sqrt{\rho}}{1-\sqrt{\rho}}\sum_{k=1}^{T-1}\eta_k,
    \end{aligned}
\end{equation*}
where the last step bounds the inner geometric series by
$\sum_{t=k+1}^{T}\rho^{(t-k)/2}\le \sqrt\rho/(1-\sqrt\rho)$. Using $\sqrt\rho/(1-\sqrt\rho)\le 2/(1-\rho)$ for $\rho\in(0,1)$ and $\eta_k=\big(\mu k+L/D\big)^{-1}$, which implies $\sum_{k=1}^{T-1}\eta_k\le \frac1\mu\log\!\big(1+\frac{T\mu D}{L}\big)$, we obtain
\begin{equation*}
\frac1n\sum_{t=1}^{T}\sum_{i=1}^n d(x_{i,t},\bar x_t)
\;\le\;\frac{2L}{\mu(1-\rho)}\log\!\Big(1+\frac{T\mu D}{L}\Big),
\end{equation*}
which provides the second bound of Lemma~\ref{lemnetworkerror}. \hfill $\square$

\subsection{Proof of Theorem \ref{thmregretboundfull}}
\label{app:thmfull}

We use two geometric lemmas. The first is a Riemannian law of cosines that two-sidedly compares squared distances. It helps with converting the gradient step into a telescoping sum and supplies the constants $C_1, C_2$. The result is due to  \cite[Corollary 2.1]{AlimisisOrvietoBecigneulLucchi2020} and \cite[Lemma~5]{ZhangSra2016}.

\begin{lem}[Trigonometric distance comparison]
\label{lemcos}
Let $a,b,c \in \mathcal{M}$ with $d(a,b)\leq D$ and $d(b,c)\leq D$, so that log maps are well-defined. Given Assumption \ref{assum:sectcurv}, we have
\begin{equation}
    \begin{aligned}
        d^2(a,c)\le&\, g_1(\kmin, d(a,b)) d^2(b,c)+d^2(a,b)\\
        &-2\langle\Log_b(a),\Log_b(c)\rangle,\\
        d^2(a,c)\ge&\, g_2(\kmax,q) d^2(b,c)+d^2(a,b)\\
        &-2\langle\Log_b(a),\Log_b(c)\rangle,
    \end{aligned}
\end{equation}
for some $q>0$, where $g_1(\cdot,\cdot)$ and $g_2(\cdot,\cdot)$ are defined in the Appendix \ref{subsecconst}. 
\end{lem}

The second result controls where positive curvature obstructs the Euclidean argument: the projection $\mathcal P_{\mathcal X}$ need not be non-expansive, so projecting the gradient update may increase the distance to $x^\ast$. The lemma bounds this excess by a gradient-norm term with constant $C_7$ \cite[Lemma~21]{WangTuHongWuShi2023}.

\begin{lem}[Projection-error bound]\label{lem:projerror}
Suppose $\mathcal{X} \subseteq M$ with diameter $D < \frac{\pi}{2\sqrt{K_{\max}}}$. Let us define the iterates $z_{i,t} = \Exp_{x_{i,t}}(-\eta_t g_t)$ and $y_{i,t+1} = P_{\mathcal{X}}(z_{i,t})$
with $\| \eta_t g_t \| \le D$. Then, it holds that
\begin{align*}
    \sum_{t=1}^T \frac{1}{2\eta_t}
    \bigl( d^2(y_{i,t+1}, x^\ast) - d^2(z_{i,t}, x^\ast) \bigr)
    \le C_7
    \sum_{t=1}^T 
    \frac{1}{2} \eta_t \|g_{t}\|^2,
\end{align*}
where $C_7$ is defined in Appendix \ref{subsecconst}.
\end{lem}

\subsubsection*{Decomposition of Regret Terms}
We now split the regret as
\begin{align}
    \Reg^{\mathrm{full}}(T)=&\frac{1}{n}\sum_{i=1}^n\sum_{t=1}^Tf_{t}(x_{i,t})-\sum_{t=1}^Tf_t(x^*)\nonumber\\
    =&\frac{1}{n}\sum_{i=1}^n\sum_{t=1}^Tf_t(x_{i,t})-\frac{1}{n}\sum_{i=1}^n\sum_{t=1}^Tf_{i,t}(x_{i,t})\label{eq:term7}\\
    &+\frac{1}{n}\sum_{i=1}^n\sum_{t=1}^Tf_{i,t}(x_{i,t})-\sum_{t=1}^Tf_t(x^*)\label{eq:term8}.
\end{align}
where \eqref{eq:term7} is the network error and \eqref{eq:term8} is the online optimization error.

\subsubsection*{Analysis of Term \eqref{eq:term7}}
Since $f_t = \frac1n \sum_{i=1}^n f_{i,t}$, we have $\frac1n \sum_{i=1}^n f_{i,t}(\bar{x}_t)=f_t(\bar{x}_t)$, so adding and subtracting $f_t(\bar{x}_t)$ leave the average unchanged as
\begin{align*}
    \frac{1}{n}\sum_{i=1}^n \left(f_{t}(x_{i,t})-f_{i,t}(x_{i,t})\right)&= \frac{1}{n}\sum_{i=1}^n \left(f_{t}(x_{i,t})-f_{t}(\bar{x}_t)\right)\\
    &+\frac{1}{n}\sum_{i=1}^n\left(f_{i,t}(\bar{x}_t)-f_{i,t}(x_{i,t})\right).
\end{align*}

By the $L$-Lipschitz property of $f_t$ and each $f_{i,t}$, both summands are at most $L·d(x_{i,t}, \bar{x}_t)$, and hence
\begin{equation*}
    \frac{1}{n}\sum_{i=1}^n \left(f_{t}(x_{i,t})-f_{i,t}(x_{i,t})\right)\le \frac{2L}{n}\sum_{i=1}^n d(x_{i,t},\bar{x}_t).
\end{equation*}
Summing over $t$ and applying Lemma \ref{lemnetworkerror}, we get
\begin{align*}
\eqref{eq:term7} \le \frac{2L}{n}\sum_{t=1}^T\sum_{i=1}^n d(x_{i,t},\bar{x}_t) \le \frac{4\,L^2}{\mu(1-\rho)}\log\!\Big(1+\frac{T\mu D}{L}\Big).   
\end{align*}

\subsubsection*{Analysis of Term \eqref{eq:term8}}
Let $z_{i,t}\coloneqq\Exp_{x_{i,t}}(-\eta_t\grad_{i,t})$ for the gradient iterate before projection. By the $\mu$-strong g-convexity of $f_{i,t}$ (Assumption \ref{assum:lossfunc}), we have
\begin{equation*}
f_{i,t}(x_{i,t})-f_{i,t}(x^*) 
\le \langle -\grad_{i,t},\Log_{x_{i,t}}(x^*)\rangle-\tfrac{\mu}{2}d^2(x_{i,t},x^*).
\end{equation*}
Since $-\grad_{i,t} = \frac{1}{\eta_t}\Log_{x_{i,t}}(z_{i,t})$, the inner product equals $\frac{1}{\eta_t}\langle \Log_{x_{i,t}}(z_{i,t}),\Log_{x_{i,t}}(x^*) \rangle $. Applying Lemma \ref{lemcos} ($a=x^\ast$, $b=x_{i,t}$,
$c=z_{i,t}$) and using the monotonicity of $g_1(K_{\min}, \cdot)$,
\begin{align*}
f_{i,t}(x_{i,t})-f_{i,t}(x^\ast)
&\le \frac{1}{2\eta_t}\big(d^2(x_{i,t},x^\ast)-d^2(z_{i,t},x^\ast)\big)\\
&+\frac{C_1}{2\eta_t}d^2(x_{i,t},z_{i,t})
-\frac{\mu}{2}d^2(x_{i,t},x^\ast).    
\end{align*}
Inserting distances based on $y_{i,t+1}=\mathcal P_{\mathcal X}(z_{i,t})$ and $x_{i,t+1}$ splits the first line into a telescoping term, a consensus-interference term, and a
projection-error term as follows
\begin{align}
f_{i,t}(x_{i,t})-f_{i,t}(x^\ast)
&\le \frac{1}{2\eta_t}\big(d^2(x_{i,t},x^\ast)-d^2(x_{i,t+1},x^\ast)\big)
\notag\\
&+\frac{1}{2\eta_t}\big(d^2(x_{i,t+1},x^\ast)-d^2(y_{i,t+1},x^\ast)\big)
\notag\\
&+\frac{1}{2\eta_t}\big(d^2(y_{i,t+1},x^\ast)-d^2(z_{i,t},x^\ast)\big)
\notag\\
&+\frac{C_1}{2\eta_t}d^2(x_{i,t},z_{i,t})
   -\frac{\mu}{2}d^2(x_{i,t},x^\ast). \label{eq:term8split}
\end{align}
The step size $\eta_t$ satisfies $\|\eta_t\grad_{i,t}\|\le\eta_1 L\le D$ for every $t$, so Lemma~\ref{lem:projerror} applies with $g_t=\grad_{i,t}$. Summing
\eqref{eq:term8split} over $t=1,\dots,T$, the first term above telescopes as
\begin{align*}
\sum_{t=1}^T\frac{1}{2\eta_t}
   &\big(d^2(x_{i,t},x^\ast)-d^2(x_{i,t+1},x^\ast)\big)\\
=&\sum_{t=2}^T\Big(\frac{1}{2\eta_t}-\frac{1}{2\eta_{t-1}}\Big)
     d^2(x_{i,t},x^\ast)\\
 &+\frac{1}{2\eta_1}d^2(x_{i,1},x^\ast)
 -\frac{1}{2\eta_T}d^2(x_{i,T+1},x^\ast).
\end{align*}
With this schedule, $\frac{1}{2\eta_t}-\frac{1}{2\eta_{t-1}}=\frac{\mu}{2}$ for all $t$, and the first line above is canceled by the $-\frac{\mu}{2}d^2(x_{i,t},x^\ast)$ term in \eqref{eq:term8split}. Using $\frac{1}{2\eta_1}-\frac{\mu}{2}=\frac{L}{2D}$ and
$d^2(x_{i,1},x^\ast)\le D^2$, we obtain
\begin{align}
\sum_{t=1}^T\big(f_{i,t}(x_{i,t})&-f_{i,t}(x^\ast)\big)
\le \frac{DL}{2}-\frac{1}{2\eta_T}d^2(x_{i,T+1},x^\ast)
\notag\\
&+\sum_{t=1}^T\frac{1}{2\eta_t}
   \big(d^2(x_{i,t+1},x^\ast)-d^2(y_{i,t+1},x^\ast)\big)
\notag\\
&+\sum_{t=1}^T(C_1+C_7)\frac{\eta_t}{2}\|\grad_{i,t}\|^2,
\label{eq:term8summed}
\end{align}
where the projection-error term has been bounded via Lemma~\ref{lem:projerror} and merged with the gradient-norm term. Since $\|\grad_{i,t}\|\le L$ and $\sum_{t}\eta_t\le\frac1\mu\log(1+\frac{T\mu D}{L})$, by discarding the non-positive term $-\frac{1}{2\eta_T}d^2(x_{i,T+1},x^\ast)$, we have
\begin{align}
&\sum_{t=1}^T\big(f_{i,t}(x_{i,t})-f_{i,t}(x^\ast)\big)
\le \frac{(C_1+C_7)L^2}{2\mu}\log\!\Big(1+\frac{T\mu D}{L}\Big)
\notag\\
&+\frac{DL}{2}+\sum_{t=1}^T\frac{1}{2\eta_t}
   \big(d^2(x_{i,t+1},x^\ast)-d^2(y_{i,t+1},x^\ast)\big).\label{eq:term8final}
\end{align}

\subsubsection*{The interference term}It remains to handle the sum term in \eqref{eq:term8final}. Since we further sum over agents in the final regret bound, we need to evaluate the following term
\begin{equation*}
\sum_{i=1}^n\sum_{t=1}^T\frac{1}{2\eta_t}
   \big(d^2(x_{i,t+1},x^\ast)-d^2(y_{i,t+1},x^\ast)\big).
\end{equation*}
For each fixed $t$ this term depends only on Step \eqref{step2}
and not on the gradient step size $\eta_t$, i.e., the factor $1/\eta_t$ is simply independent of $i$. Following the argument of \cite[proof of Thm.~IV.2]{SahinogluShahrampour2025Decentralized}, with $s=(2C_1)^{-1}C_2$ the doubly stochastic, symmetric weights make the per-$t$ sum non-positive, i.e.,
\begin{equation}
\sum_{i=1}^n\frac{1}{2\eta_t}
   \big(d^2(x_{i,t+1},x^\ast)-d^2(y_{i,t+1},x^\ast)\big)\le 0.
\label{eq:interference}
\end{equation}

\subsubsection*{Conclusion}
Combining the bounds on \eqref{eq:term7} and \eqref{eq:term8}, and using \eqref{eq:interference} to discard the interference sum in \eqref{eq:term8final}, we derive
\begin{align*}
\Reg^{\mathrm{full}}(T)
\le R_1\log\!\Big(1+\frac{T\mu D}{L}\Big)+R_2,   
\end{align*}
with $R_1:=\frac{4L^2}{\mu(1-\rho)}+\frac{(C_1+C_7)L^2}{2\mu}$ and $R_2:=\frac{DL}{2}$, both independent of $T$. This completes the proof of  Theorem~\ref{thmregretboundfull}.  \hfill $\square$

\subsection{Proof of Lemma \ref{musubconvlem}}
Let us set $x_u\coloneqq\Exp_x(u)$ and $y_u\coloneqq\Exp_y(\ptg{x}{y}(u))$, where $\ptg{x}{y}$ denotes parallel transport along the minimizing geodesic from $x$ to $y$.

\subsubsection*{Step 1 - Surrogate gap from strong convexity} Since $\ptg{x}{y}$ is a linear isometry, the change of variables $u \rightarrow \ptg{x}{y} (u)$ is measure preserving, so both surrogate integrals can be written over the common ball $\mathcal{B}_{T_{x}\mathcal{M}}(\delta)$ as follows
\begin{equation*}
    \hat{f}^{\delta}(y)-\hat{f}^{\delta}(x) = \int_{\mathcal{B}_{T_{x}\mathcal{M}}(\delta)} (f(y_u)-f(x_u)) dp(u), 
\end{equation*}
where $p$ denotes the uniform measure on $\mathcal{B}_{T_{x}\mathcal{M}}(\delta)$.
Applying $\mu$-strong g-convexity of $f$ to the pair $(x_u,y_u)$ gives
\begin{align}\label{eq:smoothSec4}
\hat{f}^{\delta}(y) - \hat{f}^{\delta}(x)
\ge& \int_{\mathcal{B}_{T_{x}\mathcal{M}}(\delta)} \langle \mathrm{grad}f(x_u),\Log_{x_u}(y_u) \rangle \ dp(u) \nonumber\\
+&\int_{\mathcal{B}_{T_{x}\mathcal{M}}(\delta)}\frac{\mu}{2}d^2(x_u,y_u) \ dp(u).
\end{align}

\subsubsection*{Step 2 - The certifying field} By the Stokes' identity for the spherical two-point estimator~\cite{AgarwalDekelXiao2010,FlaxmanKalaiMcMahan2005}, the expectation of $\esgrad$ over $u$ equals the $\delta$-ball average of the gradient of $h_x(u)\coloneqq f(\Exp_x(u))$. For any $v\in T_{x}\mathcal{M}$, we have
\begin{align}\label{eq:Stokes}
&\langle\mathbb{E}_u [\esgrad(x)],v\rangle = \nonumber\\
&\int_{\mathcal{B}_{T_{x}\mathcal{M}}(\delta)} \langle \grad{f}(x_u),d \Exp_x(u) [v] \rangle dp(u),
\end{align}
where $d\Exp_x(u):T_{x}\mathcal{M} \rightarrow T_{x_u}\mathcal{M}$ is the differential of the exponential map.

\subsubsection*{Step 3 - Subtract} Taking $v=\Log_x(y)$ in Equation \eqref{eq:Stokes} and subtracting \eqref{eq:Stokes} and $\frac{\mu}{2}d^2(x,y)$ from $\eqref{eq:smoothSec4}$, we get
\begin{align*}
    &\hat{f}^{\delta}(y) - \hat{f}^{\delta}(x) - \langle \mathbb{E}_u [\esgrad(x)], \Log_{x}(y) \rangle -\frac{\mu}{2}d^2(x,y)\\
    \ge& \int_{\mathcal{B}_{T_{x}\mathcal{M}}(\delta)} \langle \mathrm{grad}f(x_u),\Log_{x_u}(y_u) - d \Exp_x(u) [v] \rangle\ dp(u)\\
    &+ \frac{\mu}{2}\int_{\mathcal{B}_{T_{x}\mathcal{M}}(\delta)} (d^2(x_u,y_u) -d^2(x,y)) dp(u).
\end{align*}

\subsubsection*{Step 4 - Bound the two integrals} 

Applying the Cauchy-Schwarz inequality, we bound the absolute value of the integrand of the first integral as
\begin{equation*}
    \|\Log_{x_u}(y_u)- d \Exp_x(u) [\Log_x(y)]\|\cdot \|\grad f(x_u)\|.
\end{equation*}
Using $\|\grad f\|\le L$ and the bound $\|\Log_{x_u}(y_u)- d \Exp_x(u) [\Log_x(y)]\|\le \delta C_4$ from~\cite{SahinogluShahrampour2025Decentralized}, this is at most $\delta C_4 L$. Hence, the first integral is bounded below by $-\delta L C_4$. For the second integral, since $d(x,x_u)\le \delta$ and $d(y,y_u)\le \delta$, the triangle inequality gives $|d(x_u,y_u)-d(x,y)|\le 2\delta$. Since $|d(x_u,y_u)+d(x,y)|\le 2D$, we get $|(d^2(x_u,y_u) -d^2(x,y))|\le 4 \delta D$, so the integral is bounded below by $-2\mu D \delta$. We obtain
\begin{equation*}
    \begin{aligned}
        &\hat{f}^{\delta}(y) - \hat{f}^{\delta}(x) - \langle \mathbb{E}_u [\esgrad(x)], \Log_{x}(y) \rangle -\frac{\mu}{2}d^2(x,y)\\
    \ge& -\delta(L C_4+2\mu D).
    \end{aligned}
\end{equation*}

\subsection{Proof of Theorem \ref{thmbanditreg}}
Throughout, $\hat f^{\delta}_{i,t}$ and $\hat f^{\delta}_{t}=\frac1n\sum_{i=1}^n \hat f^{\delta}_{i,t}$
denote the smoothed local and global losses, $\esgrad_{i,t}$ the two-point estimator of \eqref{eq:gradestimator}, and
\begin{equation*}
x^{*}_{\tau}\;:=\;\arg\min_{x\in(1-\tau)\mathcal X}\ \sum_{t=1}^{T} f_t(x).   \end{equation*}
Recall the bandit step size $\eta_t=\bigl(\mu(t+\tfrac{dL}{\mu D})\bigr)^{-1}$ and the shrinkage factor $\tau=\delta\theta/r$. 

\subsubsection*{Step 1: Regret Decomposition}
Adding and subtracting $f_t(x_{i,t})$, $\hat f^{\delta}_t(x_{i,t})$,
$\hat f^{\delta}_t(x^{*}_{\tau})$, and $f_t(x^{*}_{\tau})$ inside the bandit regret \eqref{defregretbandit} gives $\mathbb E\!\left[\Reg^{\mathrm{bandit}}(T)\right]$ equals to
\begin{align}
&\mathbb E\,\left[\frac1n\!\sum_{i=1}^n \sum_{t=1}^T\!   \Bigl(\tfrac12\bigl(f_t(x_{i,t,1})+f_t(x_{i,t,2})\bigr)-f_t(x_{i,t})\Bigr)\right] \tag{15a}\label{15a}\\
&\quad+\mathbb E\,\left[\frac1n\!\sum_{i=1}^n \sum_{t=1}^T\!\bigl(f_t(x_{i,t})-\hat f^{\delta}_t(x_{i,t})\bigr)\right]
\;\tag{15b}\label{15b}\\
&\quad+\mathbb E\,\left[\frac1n\!\sum_{i=1}^n \sum_{t=1}^T\!\bigl(\hat f^{\delta}_t(x_{i,t})-\hat f^{\delta}_t(x^{*}_{\tau})\bigr)\right] \tag{16}\label{16}\\
&\quad+\mathbb E\,\left[\frac1n\!\sum_{i=1}^n \sum_{t=1}^T\!\bigl(\hat f^{\delta}_t(x^{*}_{\tau})- f_t(x^{*}_{\tau})\bigr)\right] \tag{15c}\label{15c}\\
&\quad+\mathbb E\,\left[\frac1n\!\sum_{i=1}^n \sum_{t=1}^T\!\bigl( f_t(x^{*}_{\tau})-f_t(x^{*})\bigr)\right]. \tag{15d}\label{15d}
\end{align}

\subsubsection*{Step 2: Upper bounds on \eqref{15a}, \eqref{15b}, \eqref{15c} and \eqref{15d}}
We show that, these terms are $O(1)$ in $T$. Each query point satisfies
$d(x_{i,t,k},x_{i,t})\le\delta$ for $k\in\{1,2\}$, so by $L$-Lipschitzness of $f_t$ we get $\eqref{15a}\le \delta T L$. Likewise,
$f_t(x)-\hat f^{\delta}_t(x)=\int_{B_{T_xM}(1)}\!\bigl(f_t(x)-f_t(\mathrm{Exp}_x(\delta u))\bigr)dp(u)\le \delta L$,
hence $\eqref{15b}\le\delta T L$ and $\eqref{15c}\le \delta T L$. 
For \eqref{15d}, $g$-convexity of $f_t$ gives
$f_t(x^{*}_{\tau})-f_t(x^{*})\le \tau\bigl(f_t(p)-f_t(x^{*})\bigr)\le \tau D L$, where $p$ is the interior reference point defined in Assumption \ref{assum:domain-sandwich} and satisfies  $\Exp_p\left((1-\tau)\Log_p(x^*)\right)\in(1-\tau)\mathcal{X}$. 
Therefore, 
\begin{equation}
\eqref{15a}+\eqref{15b}+\eqref{15c}+\eqref{15d}\;\le\;(3\delta+\tau D)\,L\,T.\tag{15}\label{15}
\end{equation}

\subsubsection*{Step 3: The main term \eqref{16}}
Insert the local smoothed losses to split \eqref{16} into a network term, an
optimization term, and a vanishing term:
\begin{align}
  \eqref{16}
&=\mathbb E\,\left[\frac1n\!\sum_{i=1}^n \sum_{t=1}^T\!\bigl(\hat f^{\delta}_t(x_{i,t})-\hat f^{\delta}_{i,t}(x_{i,t})\bigr)\right] \tag{17}\label{18}\\  
&+\mathbb E\,\left[\frac1n\!\sum_{i=1}^n \sum_{t=1}^T\!\bigl(\hat f^{\delta}_{i,t}(x_{i,t})-\hat f^{\delta}_{i,t}(x^{*}_{\tau})\bigr)\right]\tag{18}\label{19}\\
&+\mathbb E\,\left[\frac1n\!\sum_{i=1}^n \sum_{t=1}^T\!\bigl(\hat f^{\delta}_{i,t}(x^{*}_{\tau})-\hat f^{\delta}_t(x^{*}_{\tau})\bigr)\right]\tag{19}\label{20}.
\end{align}
Term \eqref{20} vanishes: $\frac1n\sum_i \hat f^{\delta}_{i,t}=\hat f^{\delta}_t$ pointwise, so
$\frac1n\sum_i\bigl(\hat f^{\delta}_{i,t}(x^{*}_{\tau})-\hat f^{\delta}_t(x^{*}_{\tau})\bigr)=0$.

\emph{Term \eqref{18} -- network error (reuse of Lemma~\ref{lemnetworkerror}).}
Firstly, we estimate the Lipschitz constant of $\hat f^{\delta}_{i,t}$. In \cite[Lemma 1]{mangoubi2018rapid} it is given that for any $w\in T_x\mathcal{M}$ satisfying $\|w\|\le\frac{\pi}{2\sqrt{\kmax}}$, we have
\begin{equation*}
    d\left(\Exp_x(w),\Exp_y(\ptg{x}{y}(w))\right)\le \cosh(\sqrt{\kmin}\|w\|)d(x,y).
\end{equation*}
Therefore, with $x_u:=\Exp_x(\delta u)$ and $y_u:=\Exp_y(\ptg{x}{y}(\delta u))$, as defined in the proof of Lemma~\ref{musubconvlem}, it follows that
\begin{equation*}
    \begin{aligned}
        \hat f^{\delta}_{i,t}(x)-\hat f^{\delta}_{i,t}(y)=&\,\int_{\mathcal{B}_{T_x\mathcal{M}}(1)}\left[f_{i,t}\left(x_u\right)-f_{i,t}\left(y_u\right)\right]dp(u)\\
        \le&\, L\cdot \int_{\mathcal{B}_{T_x\mathcal{M}}(1)} d\left(x_u,y_u\right)dp(u)\\
        \le&\, L d(x,y)\\
        &\int_{\mathcal{B}_{T_x\mathcal M}(1)} \cosh\left(\sqrt{\kmin}\delta\|u\|\right)\,dp(u) \\
        \le&\, C_8Ld(x,y),
    \end{aligned}
\end{equation*}
where $C_8$ is defined in Appendix \ref{subsecconst}. The last inequality
follows from the monotonicity of $\cosh(\cdot)$ on $[0,\infty)$ and the fact
that $\delta\|u\|\le 1$ for every $u\in \mathcal{B}_{T_x\mathcal M}(1)$. Function $\hat f^{\delta}_t$ inherits the $C_8L$-Lipschitz property from $\hat f^{\delta}_{i,t}$. Similar to the term \eqref{eq:term7} in the proof of Theorem~\ref{thmregretboundfull}, we have
\begin{equation*}
\eqref{18}\;\le\;\frac1n\sum_{i=1}^n\sum_{t=1}^T 2C_8L\,d(x_{i,t},\bar x_t).   
\end{equation*}
Lemma~\ref{lemnetworkerror} carries over almost verbatim to the bandit update \eqref{bandstep1'}: the consensus step \eqref{step2} is unchanged, and the only modification is the per-step displacement that feeds the recursion. Writing the pre-consensus chain $z_{i,t}=\mathrm{Exp}_{x_{i,t}}(-\eta_t \esgrad_{i,t})\to w_{i,t+1}=P_{\mathcal X}(z_{i,t})\to y_{i,t+1}=P_{(1-\tau)\mathcal X}(w_{i,t+1})$, the nearest-point property of the two projections, with $x_{i,t}\in(1-\tau)\mathcal X\subseteq\mathcal X$ feasible for both, gives $d(z_{i,t},w_{i,t+1})\le d(z_{i,t},x_{i,t})\le\eta_t dL$ and $d(w_{i,t+1},y_{i,t+1})\le d(w_{i,t+1},x_{i,t})\le 2\eta_t dL$ so $d(x_{i,t},y_{i,t+1})\le 4\eta_t dL$: the gradient-case displacement $\eta_t L$ is replaced by $4\eta_t dL$ (no non-expansiveness is used; on a positively curved manifold the projection onto a g-convex set may be expansive). The displacement enters the recursion linearly, enlarging only the absolute constant, not the order, so Lemma~\ref{lemnetworkerror} yields
\begin{equation}
  \eqref{18}\;\le\;\frac{16C_8d L^{2}}{(1-\rho)}\,\Big(\sumstep\Big)\le \frac{16C_8 d L^{2}}{\mu(1-\rho)}\log\left(1+\frac{T\mu D}{dL}\right).\nonumber
  \label{eq:bandit-18}
\end{equation}

\emph{Term \eqref{19} -- optimization error.}
Let $z_{i,t}$, $w_{i,t+1}$ and $y_{i,t+1}$ be defined as before. By Lemma \ref{musubconvlem}, for any fixed $t\ge1$, we have:
\begin{align*}
    &\hat{f}^{\delta}_{i,t}(x_{i,t})-\hat{f}^{\delta}_{i,t}(x^*_{\tau})\\
    \le& \langle -\mathbb{E}_u[\mathbf{g}_{i,t}^{\delta}], \Log_{x_{i,t}}(x_{\tau}^*)\rangle-\frac{\mu}{2}d^2(x_{i,t},x_{\tau}^*) +\delta\, (LC_4+2\mu D)\\
    \le& \frac{1}{2\eta_t}\mathbb{E}_u\left[d^2(x_{i,t},x_{\tau}^*)-d^2(z_{i,t},x^*_{\tau})\right]\\
    &+\frac{C_1}{2\eta_t}\cdot \mathbb{E}_u[d^2(x_{i,t},z_{i,t})]-\frac{\mu}{2}d^2(x_{i,t},x^*_{\tau})\\
    &+\delta\, (LC_4+2\mu D)\\
    =& \frac{1}{2\eta_t}\mathbb{E}_u\left[d^2(x_{i,t},x^*_{\tau})-d^2(x_{i,t+1},x^*_{\tau})\right]\\
    &+\frac{1}{2\eta_t}\mathbb{E}_u\left[d^2(x_{i,t+1},x^*_{\tau})-d^2(y_{i,t+1},x^*_{\tau})\right]\\
     &+\frac{1}{2\eta_t}\mathbb{E}_u\left[d^2(y_{i,t+1},x^*_{\tau})-d^2(w_{i,t+1},x^*_{\tau})\right]\\
    &+\frac{1}{2\eta_t}\mathbb{E}_u\left[d^2(w_{i,t+1},x^*_{\tau})-d^2(z_{i,t},x^*_{\tau})\right]\\
    &+\frac{C_1}{2\eta_t}\cdot \mathbb{E}_u[d^2(x_{i,t},z_{i,t})]-\frac{\mu}{2}d^2(x_{i,t},x^*_{\tau})\\
    &+\delta\, (LC_4+2\mu D).
\end{align*}

To apply Lemma \ref{lem:projerror}, we take $\mathbf{g}_{i,t}^{\delta}$ as the sequence $g_t$ there, and set $\eta_t = \left(\mu t + \frac{dL}{D}\right)^{-1}$ such that for each $t$, $\|\eta_t\mathbf{g}_{i,t}^{\delta}\|\le |\eta_1|\cdot dL\le D$. Then, by inequality \eqref{eq:interference} of the proof of Theorem~\ref{thmregretboundfull} with $x^{*}_{\tau}$ in place of $x^{*}$, and by Lemma \ref{lem:projerror}, we have, respectively,
\begin{equation*}
    \sum_{i=1}^n\sum_{t=1}^T\frac{1}{2\eta_t}\left(d^2(x_{i,t+1},x^*_{\tau})-d^2(y_{i,t+1},x^*_{\tau})\right)\le0,
\end{equation*}
\begin{equation*}
 \sum_{t=1}^T\frac{1}{2\eta_t}\bigl( d^2(w_{i,t+1}, x^*_{\tau}) - d^2(z_{i,t}, x^*_{\tau}) \bigr)\\
        \le C_7\sum_{t=1}^T\frac{1}{2} \eta_t \|\mathbf{g}_{i,t}^{\delta}\|^2,
\end{equation*}
and $d^2(x_{i,t},z_{i,t})=\eta_t^2\|\mathbf{g}_{i,t}^{\delta}\|^2$. Since $d(y_{i,t+1},w_{i,t+1})\le \tau D$, we also have
\begin{equation*}
    \sum_{t=1}^T\frac{1}{2\eta_t}\bigl( d^2(y_{i,t+1},x^*_{\tau})-d^2(w_{i,t+1},x^*_{\tau}) \bigr)\\
    \le \tau D^2\sum_{t=1}^T\frac{1}{\eta_t}.
\end{equation*}
Because of the term $\sum_{t=1}^T\frac{1}{\eta_t}$, both $\tau$ and $\delta$ must be $O(1/T^2)$. Set $\delta=\frac{1}{T^2}$ and keep $\eta_t = \left(\mu t + \frac{dL}{D}\right)^{-1}$. Following the same telescope summation argument as in the case of full information feedback, we have
\begin{equation*}
\eqref{19}\le R_3'\log\!\Big(1+\frac{T\mu D}{dL}\Big) + R_4',   
\end{equation*}
where $R_3'\coloneqq\frac{(C_1+C_7)(dL)^2}{2\mu}$ and $R_4'\coloneqq \frac{DdL}{2}+LC_4+2\mu D + \frac{\theta D^2}{r}(\mu+\frac{dL}{\mu D})$. Combining upper bounds on \eqref{15}, \eqref{18} and \eqref{19} we obtain 
\begin{align*}
    E\bigl[\Reg^{\mathrm{bandit}}(T)\bigr]
  \;\le\; R'_1\,\log\!\Big(1+T\frac{\mu D}{dL}\Big)\;+\;R'_2,
\end{align*}
where $R_1'\coloneqq R_3'+\frac{16C_8d L^{2}}{\mu(1-\rho)}$ and $R_2'\coloneqq R_4'+(3+\tfrac{\theta D}{r})L$. \hfill $\square$

\subsection{Constant Terms}\label{subsecconst}
In this subsection, we define the constant terms used in our paper. The first set of constants are defined as functions:
\begin{align*}
    &g_1(K,r) := \begin{cases}
    \frac{\sqrt{-K}r}{\tanh(\sqrt{-K}r)} &  \text{ if } K < 0\\
    1 &  \text{ if } K \geq 0
    \end{cases} \\
    &g_2(K,r) := \begin{cases}
    1 &  \text{ if } K\leq 0\\
    \sqrt{K}r \cot(\sqrt{K}r) &  \text{ if } K > 0
    \end{cases}
    \\
    &g_3(K,r) := \begin{cases}
    \frac{1}{6} & \text{ if } K\leq 0\\
    \frac{1}{K r^2} (\frac{r\sqrt{K}}{\sin{(r\sqrt{K})}}-1) & \text{ if } K > 0
    \end{cases}
    \\
    &g_4(K,r) := \begin{cases}
    \frac{1}{-K r^2}(\frac{\sinh(r\sqrt{-K})}{r\sqrt{-K}}-1) & \text{ if } K < 0\\
    \frac{1}{6} & \text{ if } K \geq 0
    \end{cases}
\end{align*}
Defining $\kappa:=\max\{|\kmin|,|\kmax|\}$, the second set of constants are the absolute constants: 
\begin{align*}
    &\vphantom{\frac{1}{2}}C_1 := g_1(\kmin,D)
    \vphantom{\frac{1}{2}} &C_2& := g_2(\kmax,D) \\
    &\vphantom{\frac{1}{2}} C_3 := \kappa g_3(\kmax,D) 
    \vphantom{\frac{1}{2}} &C_4& := (2C_3+C_6+3C_5)\,D^2 \\
    &\vphantom{\frac{1}{2}}C_6 :=\kappa g_4(\kmin,\delta)
    \vphantom{\frac{1}{2}} &C_7& := \max(0,-g_2(\kmax,2D))\\
    &\vphantom{\frac{1}{2}}C_8 := \cosh(\sqrt{\kmin})
\end{align*}
For $C_5$ find the definition in Lemma 6 of \cite{SunFlammarionFazel2019}.

Finally, for completeness, we present the constants in the regret bound here as well:
\begin{alignat*}{2}
R_1 &\coloneqq
    \frac{4L^2}{\mu(1-\rho)}
    +\frac{(C_1+C_7)L^2}{2\mu}
\quad
& R_2 &\coloneqq \frac{DL}{2} \\
R_1' &\coloneqq
    R_3' + \frac{16C_8d L^{2}}{\mu(1-\rho)}
& R_2' &\coloneqq
    R_4' + \left(3+\tfrac{\theta D}{r}\right)L \\
R_3' &\coloneqq
    \frac{(C_1+C_7)(dL)^2}{2\mu}
& & \\
R_4' &\coloneqq
    \mathrlap{
    \frac{DdL}{2}+LC_4+2\mu D
    +\frac{\theta D^2}{r}
    \left(\mu+\frac{dL}{\mu D}\right)}
& &
\end{alignat*}

\bibliographystyle{IEEEtran}
\bibliography{references}

\begin{IEEEbiography}[{\includegraphics[width=0.95in,height=1in,clip,keepaspectratio]{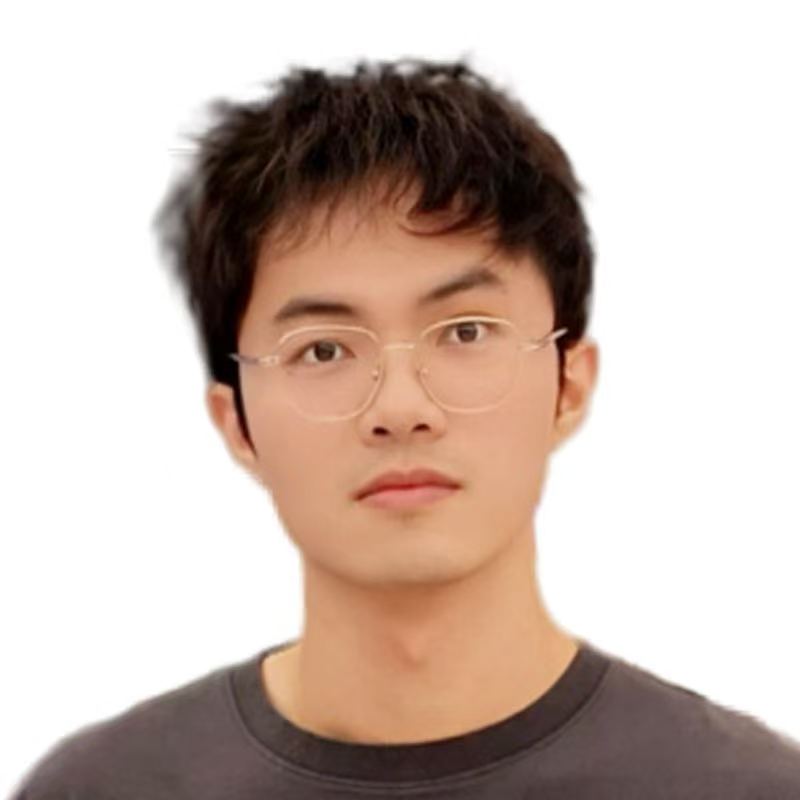}}]{Zhanyuan Cai} is currently a Ph.D. student in the Department of Mechanical and Industrial Engineering at Northeastern University. He received his B.S. degree in Mathematics and Applied Mathematics from Ocean University of China in 2022, and his M.S. degree in Mathematics from Tianjin University in 2025. His research interests include online optimization, decentralized optimization, and Riemannian optimization, with emphasis on decentralized online optimization and Riemannian geometry.
\end{IEEEbiography}

\begin{IEEEbiography}	[{\includegraphics[width=0.95in,height=1in,clip,keepaspectratio]{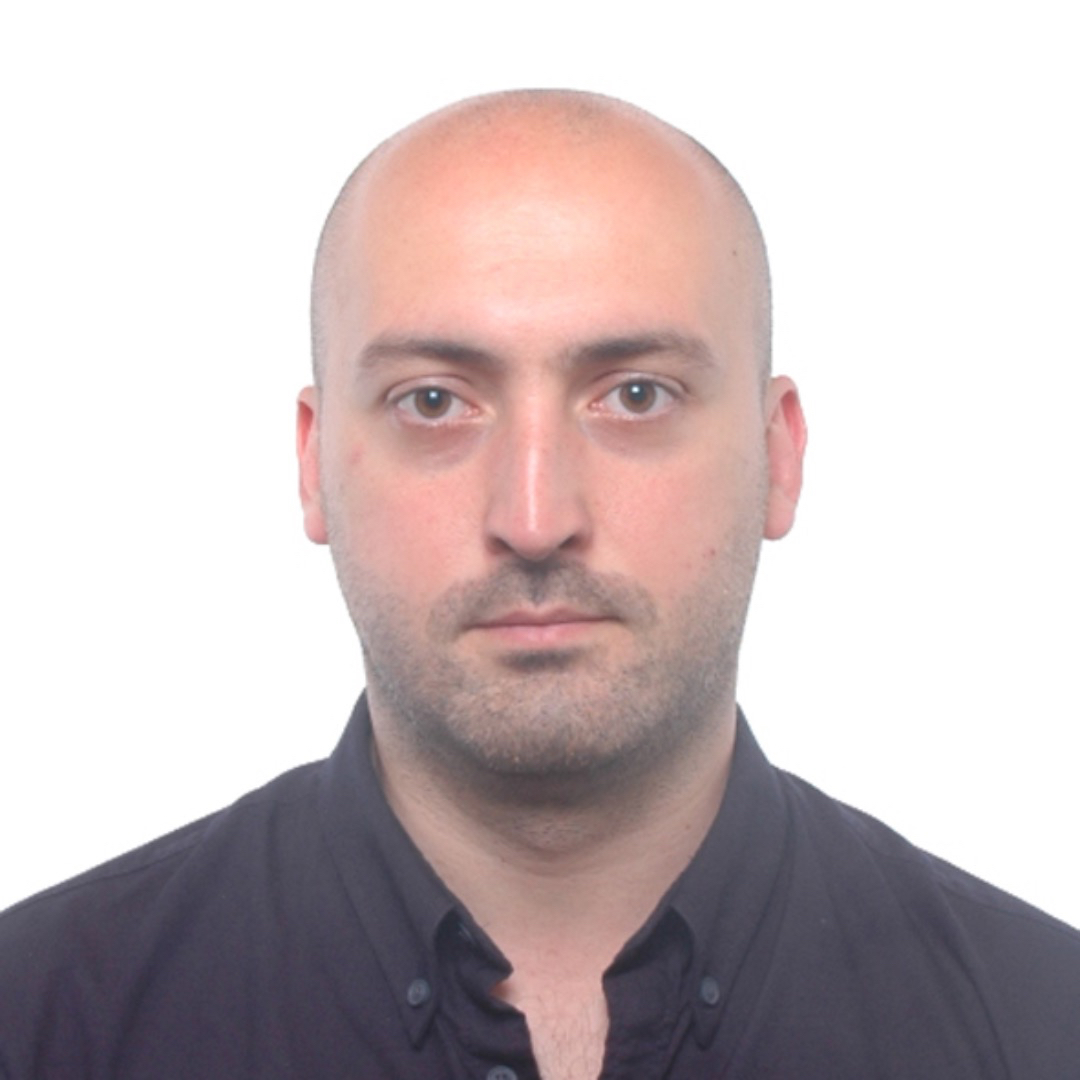}}] 
	{Emre Sahinoglu} is currently a Ph.D. candidate in the Department of Mechanical and Industrial Engineering at Northeastern University. He received the B.S. degree in Electrical and Electronics Engineering from Bilkent University, Turkey. His research interests lie at the intersection of machine learning, optimization, and multi-agent systems. In particular, his work focuses on distributed optimization, online optimization, and Riemannian optimization, with applications to learning and decision-making over networks.
\end{IEEEbiography}

\begin{IEEEbiography}
[{\includegraphics[width=1.35in,height=1.1in,clip,keepaspectratio]{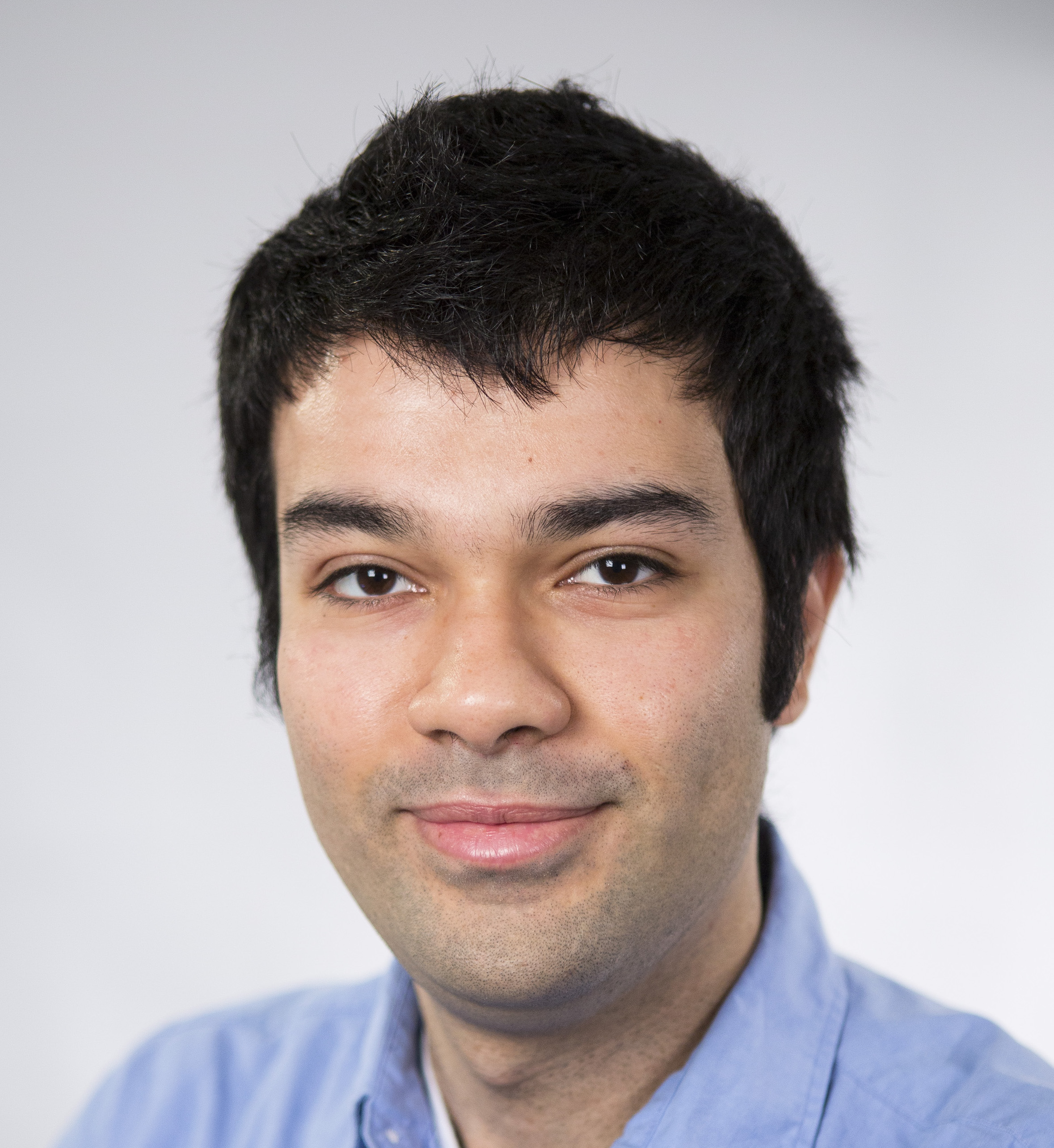}}] {Shahin Shahrampour} received the Ph.D. degree in Electrical and Systems Engineering, the M.A. degree in Statistics (The Wharton School), and the M.S.E. degree in Electrical Engineering, all from the University of Pennsylvania, in 2015, 2014, and 2012, respectively. He is currently an Associate Professor in the Department of Mechanical and Industrial Engineering at Northeastern University. His research interests include machine learning, optimization, sequential decision-making, and distributed learning, with a focus on developing computationally efficient methods for data analytics. He is a Senior Member of the IEEE.
\end{IEEEbiography}

\end{document}